\documentclass[12pt,leqno,psfig,openbib,final]{article}
\usepackage{amssymb,amsmath,graphicx,amsfonts,euscript}
\usepackage{epsfig}
\usepackage{epstopdf}
\usepackage{setspace}
\usepackage{subfigure}
\usepackage[all]{xy}
\usepackage[color]{showkeys}
\usepackage{fancyhdr}
\usepackage{hyperref}
\hypersetup{linktocpage}
\hypersetup{
    colorlinks,
    citecolor=black,
    filecolor=black,
    linkcolor=black,
    urlcolor=black
}

\newtheorem{thm}{Theorem}[section]

\newtheorem{rem}[thm]{Remark}
\newtheorem{lem}[thm]{Lemma}

\definecolor{refkey}{gray}{.75}%

\def\M{\mathcal{M}}
\def\bs{\boldsymbol}

\def\vb{\bs{v}}

\def\k{\bs{\mathrm{k}}}
\def\udiv{\mathrm{div}}

\numberwithin{equation}{section}

\begin{document}

\title{Numerical Approximation of the Inviscid 3D Primitive Equations in a Limited Domain}
\author{Qingshan Chen \thanks{Corresponding author. Mail: Department of Scientific
  Computing, Florida State University, Tallahassee, FL 32306. Email:
  {\tt qchen3@fsu.edu}}, Ming-Cheng Shiue, Roger Temam,\\ and Joseph
  Tribbia}

\date{ }
\maketitle

\thispagestyle{fancy}
\renewcommand{\headrulewidth}{0pt}
\fancyhead{ }
\fancyhead[LO]{\textsf{Article to appear in {\em M2AN Math. Model. Numer. Anal.}}} 

\begin{abstract}
   A new set of nonlocal boundary conditions are proposed
   for the higher modes of the 3D inviscid primitive equations.
   Numerical schemes using the splitting--up method are proposed
   for these modes. Numerical simulations of the full
   nonlinear primitive equations are performed on a nested
   set of domains, and the results are discussed.
\end{abstract}


\section{Introduction }\label{sec1}
When the viscosity is present, the primitive equations have been
the object of much attention, on the mathematical side. See
the original articles \cite{LTW92a, LTW92b}, and the review articles about
the mathematical theory
of the PEs with viscosity appearing in \cite{TZ04} and in an
updated form in \cite{PTZ08}; see also the articles \cite{CT07, Ko06, Ko07}.
For the physical background on primitive equations, see
e.g.~\cite{P87} or \cite{WP05}.
In the absence of viscosity, little progress has been made
on the analysis of the primitive equations since the negative
result of Oliger and Sundstr\"om \cite{OS78} showing that these
equations are not well-posed for any set of local boundary
conditions. However, the determination of suitable boundary
conditions for the primitive equations is a very important
problem for limited area models; see e.g.~a discussion
in \cite{WPT97}.

In the broader context of limited-area numerical weather prediction
modeling, the issue concerning the boundary conditions on the artificial
boundaries has been the focus of much research effort for decades.
Since the boundaries are artificial, the boundary conditions are expected
to be of a transparent, or at least, non-reflecting type. It has been pointed
out in the reference \cite{OS78} cited above, and also in \cite{EnMa77}, that
for many hyperbolic systems, truely non-reflecting boundary conditions
will have to be non-local, i.e.~involving states of the whole or part of
the time period and/or the spatial domain. Generally speaking,
non-local boundary conditions are difficult to implement in
numerical simulations. Many authors have therefore proposed
approximately non-reflecting boundary conditions. See, for example,
\cite{EnMa77}, \cite{GivNet03}, \cite{Hig86}, \cite{Mc03}, and
the references therein. These boundary conditions are also
called absorbing boundary conditions, because they are designed
to absorb the incident waves. Another approach that has been undertaken
by some authors (see e.g.~\cite{NNH04} and the references therein)
is to introduce an absorbing layer, notably
the perfectly matched layer (PML), surrounding the limited area. In
this layer, the governing equations are modified to absorb any
spurious reflections. Our approach differs from those mentioned
above in that we seek the truly non-reflecting boundary condtions,
which are suitable for the governing equations in the sense of well-posedness,
and are of a transparent type.
As we have mentioned earlier, this type of boundary conditions
will necessarily be non-local. They are valuable only if they
can be shown to be practical for implementations in limited-area
simulations. This is the main task of the current work.

Following \cite{TT03}, two of the present authors (RT and JT) and
A.~Rou\-sseau have investigated the inviscid primitive equations
in space dimension two, and
an infinite set of boundary conditions has been proposed.
Well-posedness of the corresponding linearized equations has
been established in \cite{RTT05b} and numerical simulations have
been performed in \cite{RTT07} for the linearized equations
and for the full nonlinear equations. Note that the nonoccurrence
of blow-up in the latter case supports the (yet unproved) conjecture
that the proposed nonlocal boundary conditions are also suitable
for the nonlinear PEs.

Pursuing this approach, three of the present authors (QC, RT and JT),
J.~Laminie, and A.~Rousseau considered a 2.5D model,
with three \emph{orthogonal} finite elements in the $y$-direction,
of the equations.
The well--posedness result for the linearized equations was established
in \cite{CLRTT07}, and the numerical simulations of the nonlinear
equations on a nested set of domains were discussed in
\cite{CTT08}.

The present article is related to the more theoretical ones
\cite{RTT08} and \cite{CST10}. In the first one \cite{RTT08},
the authors
obtained an infinite set of
nonlocal boundary conditions for the 3D inviscid primitive equations by
studying the stationary problem associated with the linearized
equations. In the second one \cite{CST10},
the authors
gave due treatment of the special zero
(barotropic) mode of the primitive equations and
established the well--posedness of the corresponding
linearized problem using the linear semi--group
theory. Various numerical schemes through the projection method
were also proposed, and the stability issue was studied for all of them.

In the present work we intend to discuss the numerical simulations
of the 3D \emph{nonlinear} inviscid primitive equations on a nested set of
domains.
After performing the normal mode expansions of the unknowns, we are
presented with an infinite set of 2D equations.
For the zero mode we use one of the schemes proposed
in \cite{CST10}, which is semi--implicit, and is derived by the
pressure--correction method.
For the higher modes, i.e.~the subcritical and supercritical modes,
we use the splitting--up method for the discretizations and advance
the unknowns along the $x$-- and $y$--directions in separate sub-steps.
It then seems natural to impose boundary conditions by characteristics
along the $x$-- and $y$--directions separately.
In the course of the article we recall the normal mode expansion
leading to the infinite system of 2D equations (two spatial
dimensions and time). Then we show how to discretize it in a
form suitable for the implementation of the boundary
conditions.

Two simulations are performed.
An initial simulation is carried out on a large domain with
homogeneous boundary conditions. Using the data from the
initial simulation as boundary
conditions, we perform a second simulation of the same equations
on a small interior domain. Then we compare these two results over
the interior domain.

Here the goals are twofold. On the one hand we want to numerically verify
whether the boundary conditions, proven suitable for the linearized equations,
are also suitable for the nonlinear equations. On the other hand,
we want to numerically verify the transparency property of the
proposed boundary conditions. Both goals are satisfactorily achieved.

The article is organized as follows.
In Section \ref{s2} we recall the 3D equations and their normal mode
expansion. The issues of boundary conditions and well--posedness
are also discussed. The numerical schemes are presented in Section \ref{s3}.
The settings for the numerical simulations, and the results of the
simulations are discussed in Section \ref{s4}.

\section{The model}\label{s2}
The 3D primitive equations, linearized around a uniformly stratified flow
(see \cite{RTT05b}, \cite{RTT07}, \cite{CLRTT07} and \cite{CTT08}), read
\begin{equation}\label{e1.1}
\begin{cases}
u_t+\bar U_0 u_x + \phi_x -fv+B(u,v,w;u)=0,\\
v_t+\bar U_0 v_x + \phi_y +fu + B(u,v,w;v)+f\bar {U}_0=0,\\
\psi_t + \bar U_0 \psi_x + N^2 w + B(u,v,w;\psi ) =0,\\
u_x+v_y+w_z =0,\\
\phi_z = \psi.
\end{cases}
\end{equation}
where $u$, $v$ and $w$ are the perturbation variables of the three
velocity components, $\phi$ is the perturbation variable of the pressure,
$\psi$ is the perturbation variable of the temperature; $f$ is the
Coriolis force parameter, $N$ is the Brunt--V\"ais\"al\"a (buoyancy)
frequency, assumed to be constant in the current study;
$B(u,v,w; \theta ) = u\theta_x + v\theta_y+ w\theta_z$ for $\theta = u, v $, or $\psi$. \\
\indent We will consider these equations in the domain $\M = \M'\times (-H,0), \, \\
 \M' = (0,L_1)\times (0,L_2),\, L_1, \, L_2, \, L_3 = H$ positive constants. Assuming flat bottom and the rigid lid hypothesis, we have
\begin{equation}\label{e1.1a}
w = 0 \text{ at } z = 0, \, -H.
\end{equation}
The boundary conditions for the other variables will be recalled and discussed below.
\subsection{Normal modes expansion}\label{s2.1}
Following \cite{OS78} and \cite{TT03},
we consider the normal mode expansion of the solutions
of the system \eqref{e1.1}. That is, we look for the solutions
written in the following form:
\begin{equation}
   \left\{\begin{aligned}\label{e1.2}
       &(u,v,\phi) = \sum_{n\ge 0}\mathcal{U}_n(z)(u_n,v_n,\phi_n)(x,y,t), \\
       &(w,\psi)=\sum_{n\ge 1} \mathcal {W}_n(z)(w_n,\psi_n)(x,y,t).
   \end{aligned}\right.
\end{equation}
Here $\mathcal{U}_n$ and $\mathcal{W}_n$ are solutions of the Sturm-Liouville boundary value problem
 \[
 \dfrac{d^2u}{dz^2}=-\lambda^2 u(z), \quad z \in (-H,0),
 \]
  respectively associated with the Neumann and Dirichlet boundary conditions. Therefore, we write the corresponding eigenfunctions as follows :
\begin{equation}\label{e1.3}
\left\{\begin{aligned}
&\lambda_n = \dfrac{n\pi}{H}, \\
&\mathcal{W}_n(z)=\sqrt{\dfrac{2}{H}} \sin {(\lambda_n z)}, \mathcal{U}_n(z)=\sqrt{\dfrac{2}{H}} \cos {(\lambda_n z)}, n \ge 1, \\
&\mathcal{U}_0(z)=\dfrac{1}{\sqrt{H}}.
\end{aligned}\right.
\end{equation}
We then substitute the expressions (\ref{e1.2}) into (\ref{e1.1}),
multiply each equation by $\mathcal{U}_n$ (or $\mathcal{W}_n$ for the 3rd and 5th equations),
and integrate in $z$ over the interval $(-H,0)$. We obtain the following systems:\\
\mbox{}
\\
{\it For $n=0$},
\begin{equation}\label{e1.4}
\left\{\begin{aligned}
&\dfrac{\partial u_0}{\partial t}+ \bar U_0 \dfrac{\partial u_0}{\partial x} + \dfrac{\partial \phi_0}{\partial x} -fv_0 + \int_{-H}^0 B(u,v,w;u)\, \mathcal{U}_0(z)\,dz = 0,\\
&\dfrac{\partial v_0}{\partial t}+ \bar U_0 \dfrac{\partial v_0}{\partial x} + \dfrac{\partial \phi_0}{\partial y} +fu_0 + \int_{-H}^0 B(u,v,w;v)\, \mathcal{U}_0(z)\,dz + f\bar U_0 \sqrt {H}= 0,\\
&\dfrac{\partial u_0}{\partial x}+\dfrac{\partial v_0}{\partial y} =0,\\
& \psi_0 = w_0 = 0,
\end{aligned}\right.
\end{equation}
\mbox{}
\\
{\it For $n \ge 1$,}
\begin{equation}\label{e1.5}
\left\{\begin{aligned}
&\dfrac{\partial u_n}{\partial t}+ \bar U_0 \dfrac{\partial u_n}{\partial x} + \dfrac{\partial \phi_n}{\partial x} -fv_n + \int_{-H}^0 B(u,v,w;u)\, \mathcal{U}_n(z)\, dz = 0,\\
&\dfrac{\partial v_n}{\partial t}+ \bar U_0 \dfrac{\partial v_n}{\partial x} + \dfrac{\partial \phi_n}{\partial y} +fu_n + \int_{-H}^0 B(u,v,w;v)\, \mathcal{U}_n(z)\, dz = 0,\\
&\dfrac{\partial \psi_n}{\partial t} + \bar U_0\dfrac{\partial \psi_n}{\partial x}-\dfrac{N^2}{\lambda_n}(\dfrac{\partial u_n}{\partial x}+ \dfrac{\partial v_n}{\partial y}) +\int_{-H}^0 B(u,v,w;\psi)\, \mathcal {W}_n(z)\, dz = 0.
\end{aligned}\right.
\end{equation}
\mbox{}
\\
The diagnostic unknowns $\phi_n$ and $w_n$ are given by
\begin{equation}
   \phi_n = -\dfrac{1}{\lambda_n}\psi_n,
   \label{e1.4a}
\end{equation}
and
\begin{equation}
   w_n = -\dfrac{1}{\lambda_n}(u_{nx} + v_{ny}).
   \label{e1.4b}
\end{equation}
\subsection{Boundary conditions and well--posedness issues}\label{s2.2}
With the notation $\vb=(u_0,\,v_0)^T$, the nonlinear
equations of the zero mode can be written as
\begin{equation}
   \left\{\begin{aligned}
  &\vb_t + \bar U_0\vb_x + f\k\times\vb +
  \triangledown\phi_0 + G_0 = 0,\\
  &\udiv \vb = 0.
  \end{aligned}\right.
   \label{e1.9}
\end{equation}
Here
\begin{equation}\label{e1.9a}
G_0=
\left(
  \begin{array}{c}
   \int_{-H}^0 B(u,v,w;u)\, \mathcal{U}_0(z)\,dz  \\
   \int_{-H}^0 B(u,v,w;v)\, \mathcal{U}_0(z)\,dz + f\bar U_0 \sqrt {H} \\
    \end{array}
\right),
\end{equation}
and $\triangledown$ and $\udiv$ are the 2D gradient
and divergence operators, respectively. Without considering other modes, the nonlinear equations of the zero mode become
\begin{equation}
   \left\{\begin{aligned}
  &\vb_t + \bar U_0\vb_x + f\k\times\vb +
  \triangledown\varphi +
  \dfrac{1}{\sqrt{H}}(\vb\cdot\triangledown)\vb = 0,\\
  &\udiv \vb = 0.
  \end{aligned}\right.
   \label{e1.9b}
\end{equation}
where $\varphi = \phi_0 + f\bar{U}_0 \sqrt {H} y$.\\
We supplement the system
\eqref{e1.9b} with the following boundary conditions:
\begin{equation}\label{e1.10}
\begin{cases}
u_0 = 0,\quad \text{ at } x = 0,\, L_1, \\
v_0 = 0,\quad \text{ at } x = 0, \textrm{ and } y =0,\,  L_2.
\end{cases}
\end{equation}
The well--posedness of the linearized system associated with \eqref{e1.9}, \eqref{e1.10}
has been studied in \cite{CST10}.
It is a standing conjecture that the boundary conditions
\eqref{e1.10} are also suitable for the nonlinear system
\eqref{e1.9}, at least for a certain period of time.

For the modes $n \ge 1$, we rewrite equation (\ref{e1.5}) in the matrix form as follows :
\begin{equation}\label{e1.6}
\dfrac{\partial U_n}{\partial t} + E_n\dfrac{\partial U_n}{\partial x} + F_n\dfrac{\partial U_n}{\partial y} + G_n=0.
\end{equation}
Here
\begin{equation}\label{e1.7}
U_n=
\left(
  \begin{array}{c}
    u_n \\
    v_n \\
    \psi_n \\
  \end{array}
\right),
E_n=
\left(
  \begin{array}{ccc}
    \bar U_0 &0& \dfrac{-1}{\lambda_n} \\
    0&\bar U_0& 0\\
    \dfrac{-N^2}{\lambda_n}&0&\bar U_0 \\
  \end{array}
\right),
F_n=
\left(
  \begin{array}{ccc}
    0 &0& 0 \\
    0&0& \dfrac{-1}{\lambda_n}\\
    0&\dfrac{-N^2}{\lambda_n}&0 \\
  \end{array}
\right),
\end{equation}
and
\begin{equation}\label{e1.8}
G_n=
\left(
  \begin{array}{c}
    -fv_n + \int_{-H}^0 B(u,v,w;u)\mathcal{U}_n(z)dz \\
    fu_n + \int_{-H}^0 B(u,v,w;v)\mathcal{U}_n(z)dz \\
    \int_{-H}^0 B(u,v,w;\psi)\mathcal {W}_n(z)dz \\
  \end{array}
\right),
\end{equation}

We write
\begin{equation}\label{e4.4}
\left(
  \begin{array}{c}
    \xi_n \\
    v_n \\
    \eta_n\\
  \end{array}
\right)
=
\left(
  \begin{array}{c}
   u_n-\dfrac{\psi_n}{N}  \\
    v_n \\
    u_n+\dfrac{\psi_n}{N}\\
  \end{array}
\right).
\end{equation}
and
\begin{equation}\label{e4.16}
\left(
  \begin{array}{c}
    u_n \\
    \alpha_n \\
    \beta_n\\
  \end{array}
\right)
=
\left(
  \begin{array}{c}
    u_n \\
   v_n+\dfrac{\psi_n}{N} \\
    v_n-\dfrac{\psi_n}{N}\\
  \end{array}
\right).
\end{equation}
We define $n_c$ as the positive integer satisfying the following relations:
\[
\dfrac{n_c \pi}{H} < \dfrac{N}{\bar{U}_0} < \dfrac{(n_c+1)\pi}{H}.
\]
We will not study the non generic case where $HN^{}/ \pi \bar {U}_0$ is an integer. We introduce the subcritical modes corresponding to $1 \le n < n_c$, and the supercritical modes corresponding to $n > n_c$.\\
For the {\it subcritical modes} ($n < n_c$)
we prescribe the following boundary conditions:
\begin{equation}
   \left\{\begin{aligned}
   \xi_n(0,y,t)=0,\\
   v_n(0,y,t) = 0,\\
   \eta_n(L_1,y,t) = 0.
   \label{e1.5a}
\end{aligned}\right.
\end{equation}
\begin{equation}
   \left\{\begin{aligned}
      \alpha_n(x,L_2,t)=0,\\
      \beta_n(x,0,t) =0.
   \label{e1.5b}
   \end{aligned}\right.
\end{equation}
And  for the {\it supercritical modes} ($n>n_c$) we prescribe
a slightly different set of boundary conditions:
\begin{equation}
   \left\{\begin{aligned}
   \xi_n(0,y,t)=0,\\
   v_n(0,y,t) = 0,\\
   \eta_n(0,y,t) = 0.
   \label{e1.5c}
\end{aligned}\right.
\end{equation}
\begin{equation}
   \left\{\begin{aligned}
      \alpha_n(x,L_2,t)=0,\\
      \beta_n(x,0,t) =0.
   \label{e1.5d}
   \end{aligned}\right.
\end{equation}
The well-posedness of the linearized system associated with \eqref{e1.5}
has been studied in \cite{RTT08} (see also \cite{RTT09}).
The boundary conditions \eqref{e1.5a}--\eqref{e1.5d} are
different from those proposed in \cite{RTT08} and \cite{RTT09}.
We believe that the well-posedness of the linearized system corresponding to \eqref{e1.5}
and supplemented with the foregoing boundary conditions
\eqref{e1.5a}--\eqref{e1.5d} can be established in the same way
as in \cite{RTT08} and \cite{RTT09}; this problem will be studied elsewhere.
We remark here that there are several sets of boundary conditions
which make the linearized system well-posed.\\
It is also a conjecture that the boundary conditions of \cite{RTT08} or \cite{RTT09} or the conditions \eqref{e1.5a}-\eqref{e1.5d} are suitable for the nonlinear equations for a certain time at least.

\section{The numerical schemes }\label{s3}
\subsection{Numerical scheme for the zero mode}\label{s3.1}
Due to its resemblance with the classical Navier--Stokes
equations and Euler equations,  we discretize \eqref{e1.9}
by the pressure--correction method, which is a modified
form of the classical projection method \cite{MaTe98},\cite{vK86},\cite{GMS06}.
This modified form of the projection method is known to provide a better approximation
of the pressure in the case of the Navier--Stokes equations
and we choose to use it here, instead of the initial form of the projection method \cite{Cho}, \cite{Tem69}.
The boundary conditions \eqref{e1.10} are different from
those for either the Navier--Stokes equations or the usual
Euler equations; the pressure--correction method
has to be adapted to the system \eqref{e1.9}.

We let $\Delta t = T/K$,
$\vb^k \approx \vb(x,\,y,\,k\Delta t)$, and
$\vb^{k+\frac{1}{2}}$ represents an intermediate
value between $\vb^k$ and $\vb^{k+1}$, etc.
At each step, the system is advanced in two
substeps:
\begin{equation}\label{e2.3}
\left\{\begin{aligned}
&\dfrac{\vb^{k+\frac{1}{2}} -\vb^k}{\Delta t} +
  \bar{U}_0 \vb_x^{k+\frac{1}{2}} +
  f \mathbf{k} \times \vb^k + \nabla \phi_0^k +
    G^k_0 = 0,\\
&\vb^{k+\frac{1}{2}} |_{x = 0} = 0,
\end{aligned}\right.
\end{equation}
Here
\begin{equation}\label{e2.3a}
G^k_0=
\left(
  \begin{array}{c}
   \int_{-H}^0 B(u^k,v^k,w^k;u^k)\, \mathcal{U}_0(z)\,dz  \\
   \int_{-H}^0 B(u^k,v^k,w^k;v^k)\, \mathcal{U}_0(z)\,dz + f\bar U_0 \sqrt {H} \\
    \end{array}
\right),
\end{equation}
and
\begin{equation}\label{e2.4}
\left\{\begin{aligned}
&\dfrac{\vb^{k+1}-\vb^{k+\frac{1}{2}}}{\Delta t} + \nabla (\phi_0^{k+1} - \phi_0^k) = 0,\\
&\nabla \cdot {\vb^{k+1}} = 0, \\
& \vb^{k+1} \cdot \mathbf{n} = 0.
\end{aligned}\right.
\end{equation}
where $\mathbf{n}$ is the outer normal vector on $\partial \M'$.

It has been shown in \cite{CST10} that if $\Delta t$,
$\Delta x$ and $\Delta y$ satisfy the following conditions:
\begin{equation}
   \dfrac{\Delta t}{(\Delta x^2+\Delta y^2)^2}\leq \dfrac{1}{c_1^2K_4},\qquad
   \Delta t\leq \dfrac{1}{8},
   \label{e2.5}
\end{equation}
where $c_1$ and $K_4$ are constants independent of $\Delta t$, $\Delta x$
and $\Delta y$,
then the partially implicit scheme \eqref{e2.3}--\eqref{e2.4} is stable.
For the details of the proof, and for the discussion of other related schemes,
we refer the reader to \cite{CST10}.

\subsection{Numerical scheme for the subcritical modes}\label{s3.2}
In this subsection we use the splitting method \cite{Mar71},\cite{Yan71},\cite{Tem69} to
discretize \eqref{e1.6}, in the case of the subcritical modes.
The supercritical case will be discussed in the next subsection.

We have seen that the domain under consideration is $(0,L_1)\times(0,L_2)\times(0,T)$.
We let $I,J$ and $K$ denote the numbers of grid pints in the $x$--direction,
the $y$--direction and in time, and let
$\Delta x = L_1/I, \Delta y = L_2/J$, and $\Delta t = T/K$, denote the corresponding mesh.
We denote the discrete grid points by $(x_i,\,y_j,\,t_k) = (i\Delta x,\,
j\Delta y,\, k\Delta t)$.
We let $U_n^k$ be the semi--discrete approximate value of $U_n$ at time
$t_k$, and $U^k_{n,i,j}$ the fully discrete approximate value of $U_n$
at $(x_i,\,y_j,\,t_k)$.

For each time step, two substeps are involved.
The first substep
is meant to advance the unknowns using only the advective terms in the
$x-$direction and the zero order terms (the Coriolis force), that is, to
solve the following semi-discrete equation:
\begin{equation}\label{e4.1}
   \dfrac{U_n^{k+\frac{1}{2}} - U_n^k}{\Delta t} +
   E_n\dfrac{\partial U_n^{k+\frac{1}{2}}}{\partial x}  + G_n = 0.
\end{equation}
The second substep is meant to advance the unknowns using only the advective
terms in the $y-$direction, that is, to solve the following
semi--discrete equation:
\begin{equation}\label{e4.2}
   \dfrac{U_n^{k+1} - U_n^{k+\frac{1}{2}}}{\Delta t} +
   F_n\dfrac{\partial U_n^{k+1}}{\partial y}  =0.
\end{equation}

We now discuss the full discretization of the equations \eqref{e4.1}
and \eqref{e4.2}. We first apply the change of variables
\eqref{e4.4} to \eqref{e4.1}, and obtain
\begin{equation}\label{e4.5}
\left\{\begin{aligned}
   \dfrac{\xi_n^{k+\frac{1}{2}} - \xi_n^k}{\Delta t}+
   (\bar U_0+\dfrac{N}{\lambda_n})
   \dfrac{\partial \xi_n^{k+\frac{1}{2}}}{\partial x}
= fv_n^k - \int_{-H}^0 B(u^k,v^k,w^k;u^k)\mathcal{U}_n(z)dz &\\
+
\dfrac{1}{N}\int_{-H}^0 B(u^k,v^k,w^k;\psi^k)\mathcal{W}_n(z)dz,& \\
\dfrac{v_n^{k+\frac{1}{2}} - v_n^k}{\Delta t} +
\bar U_0 \dfrac{\partial v_n^{k+\frac{1}{2}}}{\partial x} =
-f\dfrac{\eta_n^k+\xi_n^k}{2} -
\int_{-H}^0 B(u^k,v^k,w^k;v^k)\mathcal{U}_n(z)dz,&\\
\dfrac{\eta_n^{k+\frac{1}{2}} - \eta_n^k}{\Delta t}+
(\bar U_0-\dfrac{N}{\lambda_n})\dfrac{\partial \eta_n^{k+\frac{1}{2}}}{\partial x}
= fv_n^k- \int_{-H}^0 B(u^k,v^k,w^k;u^k)\mathcal{U}_n(z)dz &\\
 -
\dfrac{1}{N}\int_{-H}^0 B(u^k,v^k,w^k;\psi^k)\mathcal{W}_n(z)dz,&
\end{aligned}\right.
\end{equation}
We recall that, for the subcritical modes, $\bar {U}_0 - N/\lambda_n < 0$.
The up-wind method applied to \eqref{e4.5} yields the following fully discrete
scheme:
\begin{equation}\label{e4.6}
\begin{cases}
\dfrac{\xi_{n,i,j}^{k+\frac{1}{2}}-\xi_{n,i,j}^k}{\Delta t} + (\bar U_0+\dfrac{N}{\lambda_n}) \dfrac{\xi_{n,i,j}^{k+\frac{1}{2}}-\xi_{n,i-1,j}^{k+\frac{1}{2}}}{\Delta x} =S_{n,i,j}^{k,1}, i=2,\cdots, I+1,  \\
\dfrac{v_{n,i,j}^{k+\frac{1}{2}}-v_{n,i,j}^k}{\Delta t} + \bar U_0 \dfrac{v_{n,i,j}^{k+\frac{1}{2}}-v_{n,i-1,j}^{k+\frac{1}{2}}}{\Delta x} =S_{n,i,j}^{k,2}, i =2,\cdots, I+1,\\
\dfrac{\eta_{n,i,j}^{k+\frac{1}{2}}-\eta_{n,i,j}^k}{\Delta t} + (\bar U_0-\dfrac{N}{\lambda_n}) \dfrac{\eta_{n,i+1,j}^{k+\frac{1}{2}}-\eta_{n,i,j}^{k+\frac{1}{2}}}{\Delta x} =S_{n,i,j}^{k,3},i=1,\cdots,I,\\
\hspace{5cm} \text{ and } j = 1, \cdots, J+1 \text{ in all cases},
\end{cases}
\end{equation}

where
\begin{equation}\label{e4.7}
\left\{\begin{aligned}
S_{n,i,j}^{k,1} =\xi_{n,i,j}^k+fv_{n,i,j}^k- \int_{-H}^0 B(u_{i,j}^k,v_{i,j}^k,w_{i,j}^k;u_{i,j}^k)\mathcal{U}_n(z)dz& \\
 +\dfrac{1}{N}\int_{-H}^0 B(u_{i,j}^k,v_{i,j}^k,w_{i,j}^k;\psi_{i,j}^k)\mathcal {W}_n(z)dz,& \\
S_{n,i,j}^{k,2} = v_{n,i,j}^k-f\dfrac{\eta_{n,i,j}^k+\xi_{n,i,j}^k}{2} - \int_{-H}^0 B(u_{i,j}^k,v_{i,j}^k,w_{i,j}^k;v_{i,j}^k)\mathcal{U}_n(z)dz,\\
S_{n,i,j}^{k,3} =\eta_{n,i,j}^k+fv_{n,i,j}^k- \int_{-H}^0 B(u_{i,j}^k,v_{i,j}^k,w_{i,j}^k;u_{i,j}^k)\mathcal{U}_n(z)dz& \\
-\dfrac{1}{N}\int_{-H}^0 B(u_{i,j}^k,v_{i,j}^k,w_{i,j}^k;\psi_{i,j}^k)\mathcal {W}_n(z)dz,&
\end{aligned}\right.
\end{equation}
The boundary conditions for
$\xi_n^{k+\frac{1}{2}}$,
$v_n^{k+\frac{1}{2}}$ and
$\eta_n^{k+\frac{1}{2}}$ are,
for $0\leq j \leq J$,
\begin{equation}
   \xi^{k+\frac{1}{2}}_{n,0,j} = 0,\quad
      v^{k+\frac{1}{2}}_{n,0,j} = 0,\quad
      \eta^{k+\frac{1}{2}}_{n,I,j} = 0.
   \label{e4.7a}
\end{equation}
We then apply the change of variables \eqref{e4.16} to \eqref{e4.2} and obtain
\begin{equation}\label{e4.17}
\left\{\begin{aligned}
   &\dfrac{u_n^{k+1} - u_n^{k+\frac{1}{2}}}{\Delta t} =0,\\
   &\dfrac{\alpha_n^{k+1}-\alpha_n^{k+\frac{1}{2}}}{\Delta t}-
   \dfrac{N}{\lambda_n}\dfrac{\partial \alpha_n^{k+1}}{\partial y} =0, \\
   &\dfrac{\beta_n^{k+1} - \beta_n^{k+\frac{1}{2}}}{\Delta t}+
   \dfrac{N}{\lambda_n}\dfrac{\partial \beta_n^{k+1}}{\partial y} =0.
\end{aligned}\right.
\end{equation}
Applying the up--wind method to the system \eqref{e4.17} yields
\begin{equation}\label{e4.18}
\left\{\begin{aligned}
&\dfrac{u_{n,i,j}^{k+1}-u_{n,i,j}^{k+\frac{1}{2}}}{\Delta t} =0, \\
&\dfrac{\alpha_{n,i,j}^{k+1}-\alpha_{n,i,j}^{k+\frac{1}{2}}}{\Delta t} -\dfrac{N}{\lambda_n}\dfrac{\alpha_{n,i,j+1}^{k+1}-\alpha_{n,i,j}^{k+1}}{\Delta y} = 0,\\
&\dfrac{\beta_{n,i,j}^{k+1}-\beta_{n,i,j}^{k+\frac{1}{2}}}{\Delta t} +\dfrac{N}{\lambda_n}\dfrac{\beta_{n,i,j}^{k+1}-\beta_{n,i,j-1}^{k+1}}{\Delta y} = 0
\end{aligned}\right.
\end{equation}
The boundary conditions for $\alpha_n^{k+1}$, $\beta_n^{k+1}$ are
\begin{equation}
   \left\{\begin{aligned}
      \alpha^{k+1}_{n,I,j}=0,\qquad\textrm{ for } 0\leq j\leq J,\\
      \beta^{k+1}_{n,i,0}=0,\qquad\textrm{ for } 0\leq i\leq I.
      \label{e4.17a}
   \end{aligned}\right.
\end{equation}
We remark here that $u_n^{k+1}$ does not need any boundary conditions.

\subsection{Numerical scheme for the supercritical modes}\label{s3.3}
The fully discrete numerical schemes for the supercritical modes can
be derived by the same approach presented in the previous subsection. The
results for the supercritical modes are similar to those for the subcritical
modes, and are simpler because all the eigenvalues of the coefficient
matrix $E_n$ are positive. We shall omit the intermediate details, and present
the numerical schemes directly. Only the differences with those for the subcritical
modes will be pointed out.

As for the subcritical modes, the numerical schemes for the supercritical
modes also involve two substeps. The first substep consists of the following
scheme:
\begin{equation}\label{e4.19}
\begin{cases}
\dfrac{\xi_{n,i,j}^{k+\frac{1}{2}}-\xi_{n,i,j}^k}{\Delta t} + (\bar U_0+\dfrac{N}{\lambda_n}) \dfrac{\xi_{n,i,j}^{k+\frac{1}{2}}-\xi_{n,i-1,j}^{k+\frac{1}{2}}}{\Delta x} =S_{n,i,j}^{k,1}, i=2,\cdots,I+1,\\
\dfrac{v_{n,i,j}^{k+\frac{1}{2}}-v_{n,i,j}^k}{\Delta t} + \bar U_0 \dfrac{v_{n,i,j}^{k+\frac{1}{2}}-v_{n,i-1,j}^{k+\frac{1}{2}}}{\Delta x} =S_{n,i,j}^{k,2}, i=2,\cdots,I+1,\\
\dfrac{\eta_{n,i,j}^{k+\frac{1}{2}}-\eta_{n,i,j}^k}{\Delta t} + (\bar U_0-\dfrac{N}{\lambda_n}) \dfrac{\eta_{n,i,j}^{k+\frac{1}{2}}-\eta_{n,i-1,j}^{k+\frac{1}{2}}}{\Delta x} =S_{n,i,j}^{k,3},i=1,\cdots,I,\\
\hspace{5cm} \text{ and } i=1,\cdots, J+1 \text{ in all cases },
\end{cases}
\end{equation}
Here, $S_{n,i,j}^{k,1},\, S_{n,i,j}^{k,2},$ and $S_{n,i,j}^{k,3}$ are defined in \ref{e4.7}.
We note here that $\partial \eta_n^{k+\frac{1}{2}}/\partial x$ is discretized differently in $\eqref{e4.6}$ and in \eqref{e4.19},
due to the fact that $\bar {U}_0 - 1/\lambda_n$ has different signs in the sub-- and super--critical modes.
The boundary conditions for
$\xi_n^{k+\frac{1}{2}}$,
$v_n^{k+\frac{1}{2}}$ and
$\eta_n^{k+\frac{1}{2}}$ are,
for $0\leq j \leq J$,
\begin{equation}
      \xi^{k+\frac{1}{2}}_{n,0,j} = 0,
      \quad v^{k+\frac{1}{2}}_{n,0,j} = 0,
      \quad \eta^{k+\frac{1}{2}}_{n,0,j} = 0.
   \label{e4.21}
\end{equation}
The second substep consists of the following scheme:
\begin{equation}\label{e4.22}
\left\{\begin{aligned}
&\dfrac{u_{n,i,j}^{k+1}-u_{n,i,j}^{k+\frac{1}{2}}}{\Delta t} =0, \\
&\dfrac{\alpha_{n,i,j}^{k+1}-\alpha_{n,i,j}^{k+\frac{1}{2}}}{\Delta t} -\dfrac{N}{\lambda_n}\dfrac{\alpha_{n,i,j+1}^{k+1}-\alpha_{n,i,j}^{k+1}}{\Delta y} = 0,\\
&\dfrac{\beta_{n,i,j}^{k+1}-\beta_{n,i,j}^{k+\frac{1}{2}}}{\Delta t} +\dfrac{N}{\lambda_n}\dfrac{\beta_{n,i,j}^{k+1}-\beta_{n,i,j-1}^{k+1}}{\Delta y} = 0
\end{aligned}\right.
\end{equation}
The boundary conditions for $\alpha_n^{k+1}$, $\beta_n^{k+1}$ are
\begin{equation}
   \left\{\begin{aligned}
      \alpha^{k+1}_{n,I,j}=0,\qquad\textrm{ for } 0\leq j\leq J,\\
      \beta^{k+1}_{n,i,0}=0,\qquad\textrm{ for } 0\leq i\leq I.
      \label{e4.23}
   \end{aligned}\right.
\end{equation}

\subsection{Treatment of the integral of the nonlinear term }\label{s3.4}
In this section, we will deal with the integral of the nonlinear term. There are five kinds of integrals to be considered. We have the following lemma.
\begin{lem}\label{l6.1}
Assume that $u, v, \phi, w$, and $\psi$ have the expressions (\ref{e1.2}) and
$\mathcal{U}_n$ and $\mathcal{W}_n$ for $n\geq 0$ are as defined in \eqref{e1.3}. Then

\begin{equation}\label{e6.1}
\int_{-H}^0 B(u,v,w;u)\mathcal{U}_0(z)dz = \dfrac{1}{\sqrt{H}} \sum_{m \ge 0} (u_m \dfrac{\partial u_m}{\partial x} + v_m\dfrac{\partial u_m}{\partial y}) -\dfrac{1}{\sqrt {H}}\sum_{m \ge 1} \lambda_m w_m u_m.
\end{equation}

\begin{equation}\label{e6.2}
\int_{-H}^0 B(u,v,w;v)\mathcal{U}_0(z)dz =\dfrac{1}{\sqrt{H}} \sum_{m \ge 0} (u_m \dfrac{\partial v_m}{\partial x} + v_m\dfrac{\partial v_m}{\partial y}) -\dfrac{1}{\sqrt {H}}\sum_{m \ge 1} \lambda_m w_m v_m.
\end{equation}

\begin{equation}\label{e6.3}
\begin{aligned}
\int_{-H}^0 B(u,v,w;u)\mathcal{U}_n(z)dz =
\dfrac{1}{\sqrt{2H}} \sum_{m \ge 0}^n ( u_{n-m}\dfrac{\partial u_m}{\partial x} + v_{n-m}\dfrac{\partial u_m}{\partial y})  \\
\dfrac{1}{\sqrt{2H}} \sum_{m =n}^{\infty} ( u_{m-n}\dfrac{\partial u_m}{\partial x} + v_{m-n}\dfrac{\partial u_m}{\partial y})
+
\dfrac{1}{\sqrt{2H}} \sum_{m \ge 0}^{\infty} ( u_{m+n}\dfrac{\partial u_m}{\partial x} + v_{m+n}\dfrac{\partial u_m}{\partial y}) \\
-
\dfrac{1}{\sqrt{2H}} \sum_{m \ge n+1}^{\infty} \lambda_m w_{m-n} u_m
-
\dfrac{1}{\sqrt{2H}} \sum_{m \ge 1}^{\infty} \lambda_m w_{m+n} u_m
+
\dfrac{1}{\sqrt{2H}} \sum_{m \ge 1}^n \lambda_m w_{n-m} u_m
\end{aligned}
\end{equation}

\begin{equation}\label{e6.4}
\begin{aligned}
\int_{-H}^0 B(u,v,w;v)\mathcal{U}_n(z)dz =
\dfrac{1}{\sqrt{2H}} \sum_{m \ge 0}^n ( u_{n-m}\dfrac{\partial v_m}{\partial x} + v_{n-m}\dfrac{\partial v_m}{\partial y}) \\
\dfrac{1}{\sqrt{2H}} \sum_{m =n}^{\infty} ( u_{m-n}\dfrac{\partial v_m}{\partial x} + v_{m-n}\dfrac{\partial v_m}{\partial y})
+
\dfrac{1}{\sqrt{2H}} \sum_{m \ge 0}^{\infty} ( u_{m+n}\dfrac{\partial v_m}{\partial x} + v_{m+n}\dfrac{\partial v_m}{\partial y}) \\
-
\dfrac{1}{\sqrt{2H}} \sum_{m \ge n+1}^{\infty} \lambda_m w_{m-n} v_m
-
\dfrac{1}{\sqrt{2H}} \sum_{m \ge 1}^{\infty} \lambda_m w_{m+n} v_m
+
\dfrac{1}{\sqrt{2H}} \sum_{m \ge 1}^n \lambda_m w_{n-m} v_m
\end{aligned}
\end{equation}
\begin{equation}\label{e6.5}
\begin{aligned}
\int_{-H}^0 B(u,v,w;\psi)\mathcal{W}_n(z)dz =
\dfrac{1}{\sqrt{2H}} \sum_{m=n}^{\infty} ( u_{m-n}\dfrac{\partial {\psi}_m}{\partial x} + v_{m-n}\dfrac{\partial {\psi}_m}{\partial y})  \\
+
\dfrac{1}{\sqrt{2H} }\sum_{m=1}^n ( u_{n-m}\dfrac{\partial {\psi}_m}{\partial x} + v_{n-m}\dfrac{\partial {\psi}_m}{\partial y})
-
\dfrac{1}{\sqrt{2H}} \sum_{m=1}^{\infty} ( u_{n+m} \dfrac{\partial {\psi}_m}{\partial x} + v_{n+m}\dfrac{\partial {\psi}_m}{\partial y}) \\
-
\dfrac{1}{\sqrt{2H}} \sum_{m \ge 1}^n \lambda_m w_{n-m}\psi_m
-
\dfrac{1}{\sqrt{2H}} \sum_{m \ge 1}^{\infty}  \lambda_m w_{n+m}\psi_m
+
\dfrac{1}{\sqrt{2H}} \sum_{m \ge n}^{\infty}  \lambda_m w_{m-n}\psi_m
\end{aligned}
\end{equation}
\end{lem}
Lemma \ref{l6.1} can be verified by direct calculations.

\begin{rem}\label{r6.1}
In large-scale GFD simulations, in which a large number of modes are involved, the preceding convolution products would be too costly in terms of CPU time to be appropriate. To avoid them, it is then necessary to transform the Fourier coefficients $u_n,$ etc., back into the physical space, compute the nonlinear products in the physical space, and calculate the integrals on the left side of (\ref{e6.1})-(\ref{e6.5}). In our study, only a small number ($\le 10$) of modes are considered, and thus the formulas (\ref{e6.1})-(\ref{e6.5}) are appropriate and sufficient.
\end{rem}

\section{Numerical simulations in a nested environment}\label{s4}
\begin{figure}
\begin{center}
\epsfig{file= 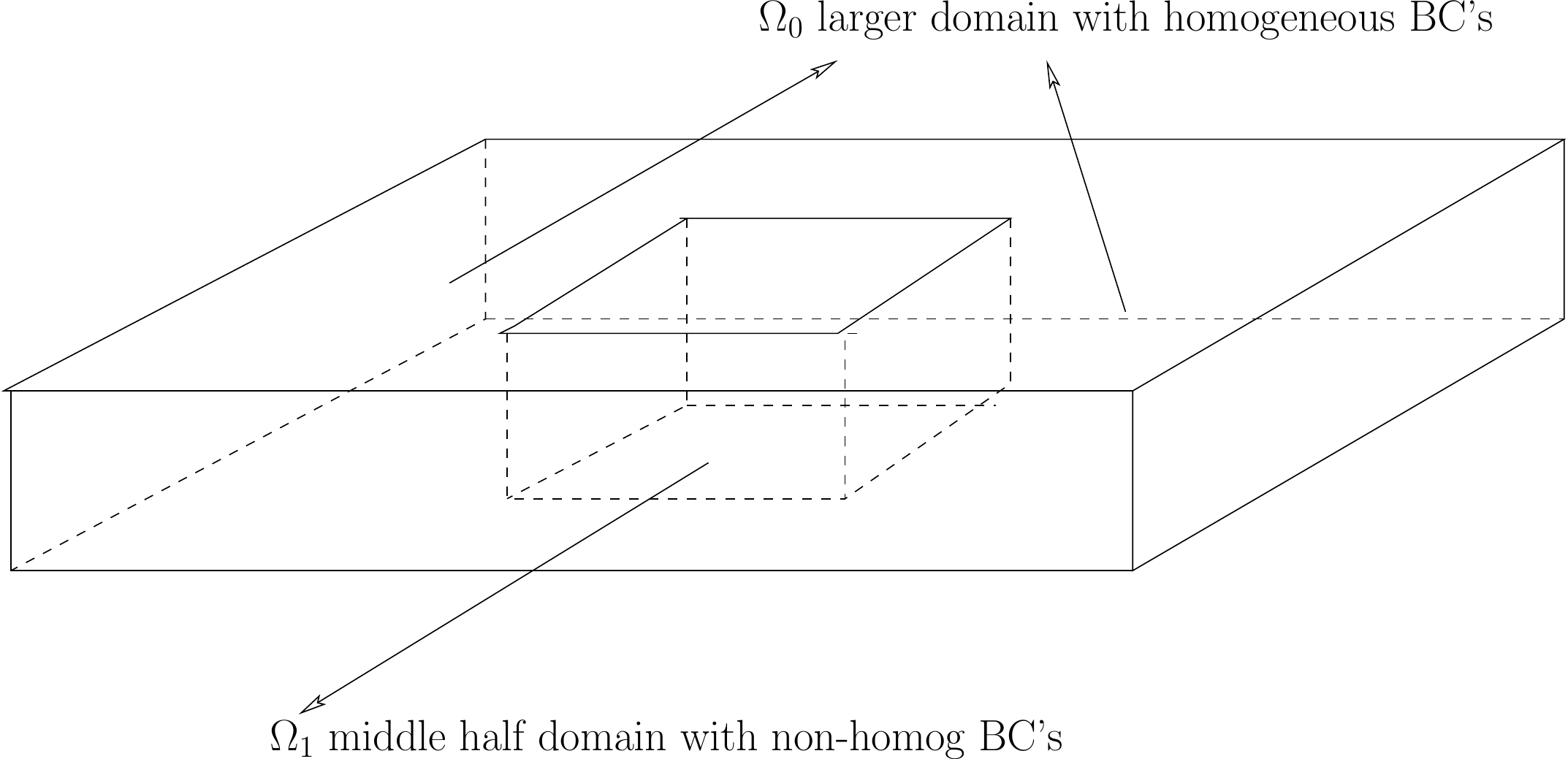, height=5cm, width=8cm}
\caption{The larger domain $\M$ and the middle half domain $\M_1$}
\label{f1.1}
\end{center}
\end{figure}

Two different simulations are performed. The first one is carried out on the larger domain $\M = (0,L_1)\times (0,L_2)\times(-H,0)$ (see Figure \ref{f1.1}), and a set of homogeneous boundary conditions prescribed at $ (x,y) \in \partial \,  \M'$, where $\M' = (0,L_1)\times (0,L_2)$. The simulations will be described and the results will be presented in details in Section \ref{s4.1}. The data obtained through this simulation will provide the nonhomogeneous boundary conditions for the second simulation on the middle half domain, denoted by $\M_1=(L_1/4,3L_1/4)\times (L_2/4,3L_2/4)\times(-H,0)$(see also Figure \ref{f1.1}), of $\M$. This simulation will be described and the numerical results will be presented in detail in Section \ref{s4.2}.

In Section \ref{s4.3}, the numerical results from these two simulations are then compared,
and the coincidence of the numerical results demonstrates the transparent properties
of the proposed boundary conditions, and supports the conjecture of their suitability for the nonlinear equations.

The physical parameters that we used in the simulations are
the following ones: $L_1 = 1000$km, $L_2 = 500$km, $H = 10$km.
We take the constant reference velocity $\bar{U}_0 = 20$ m/s,
the Coriolis parameter $f = 10^{-4}$, and the Brunt-V\"ais\"al\"a
(buoyancy) frequency $N = 10^{-2}$. The final time for
the simulations is $T = 5 \times 10^4$s, and we take
1600 time steps. In the vertical direction we take 40 segments.
In the computations, we will deal with $N_{\text{max}} = 5$ (the number of modes), which is sufficient from the physical point of view.

\subsection{Simulation on the larger domain}\label{s4.1}
In the simulation, the initial conditions are given for these scalar functions:
\begin{equation}\label{e7.1}
\begin{cases}
u(x,y,z,0) = \dfrac{x}{L_1}\dfrac{2\pi}{L_2} \sin{(\dfrac{2\pi x}{L_1})} \cos{(\dfrac{2\pi y}{L_2})} + \sin{(\dfrac{4\pi x}{L_1})} \cos{(\dfrac{4\pi y}{L_2})} \cos{(\dfrac{\pi z}{H})},\\
v(x,y,z,0)=\dfrac{-1}{L_1}\left (\sin{(\dfrac{2\pi x}{L_1})}+ \dfrac{2\pi x}{L_1}\cos{(\dfrac{2\pi x}{L_1})}\right )\sin{(\dfrac{2\pi y}{L_2})}\\
\hspace{2.1 cm}  +\dfrac{L_2}{L_1}\left ( \sin^2{(\dfrac{4\pi x}{L_1})}+ \sin{(\dfrac{4\pi x}{L_1})} \sin{(\dfrac{4\pi y}{L_2})} \cos{(\dfrac{\pi z}{H})}\right ),\\
w(x,y,z,0)=\dfrac{-4H}{L_1}(\sin{(\dfrac{4\pi x}{L_1})}+\cos{(\dfrac{4\pi x}{L_1})})\cos{(\dfrac{4\pi y}{L_2})} \sin{(\dfrac{\pi z}{H})},\\
\phi(x,y,z,0)=\bar{U}_0 \sin{(\dfrac{2\pi x}{L_1})} \sin{(\dfrac{2\pi y}{L_2})}(\cos{(\dfrac{\pi z}{H})} - \cos{(\dfrac{2\pi z}{H})}),\\
\psi(x,y,z,0)=\dfrac{\pi \bar{U}_0}{H}\sin{(\dfrac{2\pi x}{L_1})} \sin{(\dfrac{2\pi y}{L_2})}(2\sin{(\dfrac{2\pi z}{H})} - \sin{(\dfrac{\pi z}{H})}).
\end{cases}
\end{equation}
We note here that these initial functions $u, v, w, \phi$, and $\psi$ satisfy the homogeneous boundary conditions for each mode $ n \ge 0$. Specifically, for the {\it zeroth mode}, i.e. when $n = 0$,
\begin{equation}\label{e7.2}
\begin{cases}
u_0(0, y, t) = 0, \quad u_0(L_1,y,t) = 0,\\
v_0(0,y,t) = 0,\quad v_0(x,0,t) = 0, \quad v_0(x,L_2,t) = 0;
\end{cases}
\end{equation}
for the {\it subcritical modes}, i.e. when $1 \le n < n_c$,
\begin{equation}\label{e7.3}
\begin{cases}
\xi_n(0,y,t) = 0,\quad v_n(0,y,t) = 0,\quad \eta_n(L_1,y,t) = 0,\\
\alpha_n(x,L_2,t) = 0,\quad \beta_n(x,0,t) =0;
\end{cases}
\end{equation}
and for the {\it supercritical modes}, i.e. when $ n > n_c$,
\begin{equation}\label{e7.4}
\begin{cases}
\xi_n(0,y,t) = 0,\quad v_n(0,y,t) = 0,\quad \eta_n(0,y,t) = 0,\\
\alpha_n(x,L_2,t) = 0,\quad \beta_n(x,0,t) =0.
\end{cases}
\end{equation}

In this simulation, we take 400 segments in the $x$-direction,
and 200 segments in the $y$ direction. When restricted to
the middle half domain, the functions \eqref{e7.1} also
provide the \emph{initial conditions} for the
simulations on the middle half domain.

\begin{figure}
\begin{center}
\includegraphics[scale=0.8]{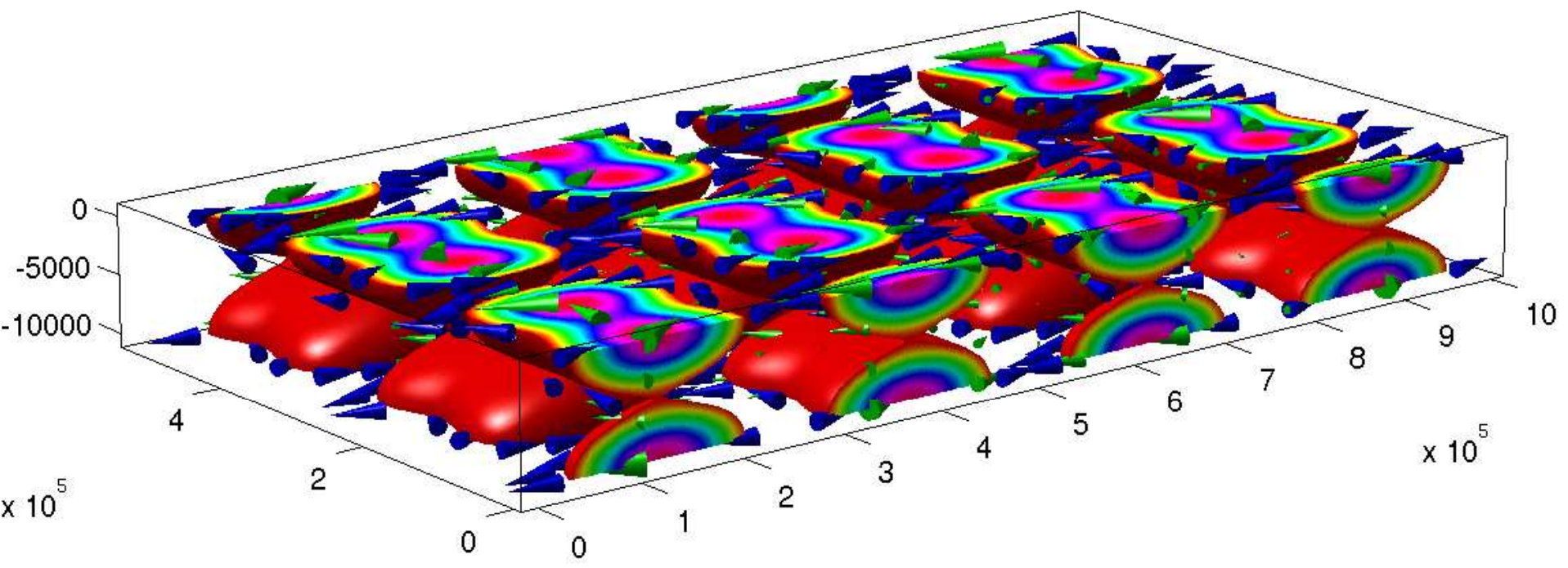}
\caption{The initial state of velocity field in the larger domain $\M$}
\label{f1.2}
\end{center}
\end{figure}

\begin{figure}
\begin{center}
\includegraphics[scale=0.9]{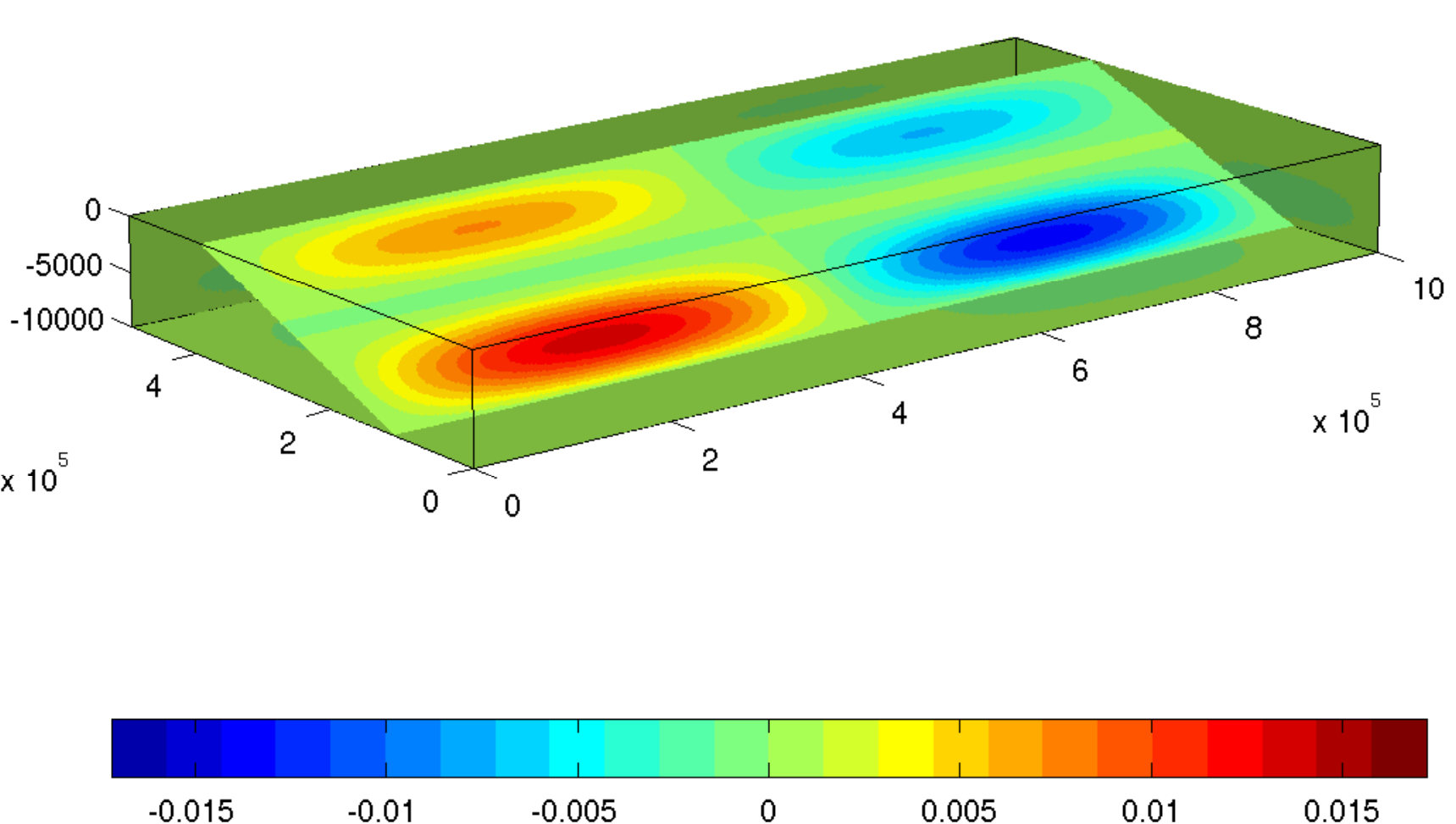}
\caption{The initial state of $\psi$ in the larger domain $\M$}
\label{f1.3}
\end{center}
\end{figure}

\begin{figure}
\begin{center}
\includegraphics[scale=0.9]{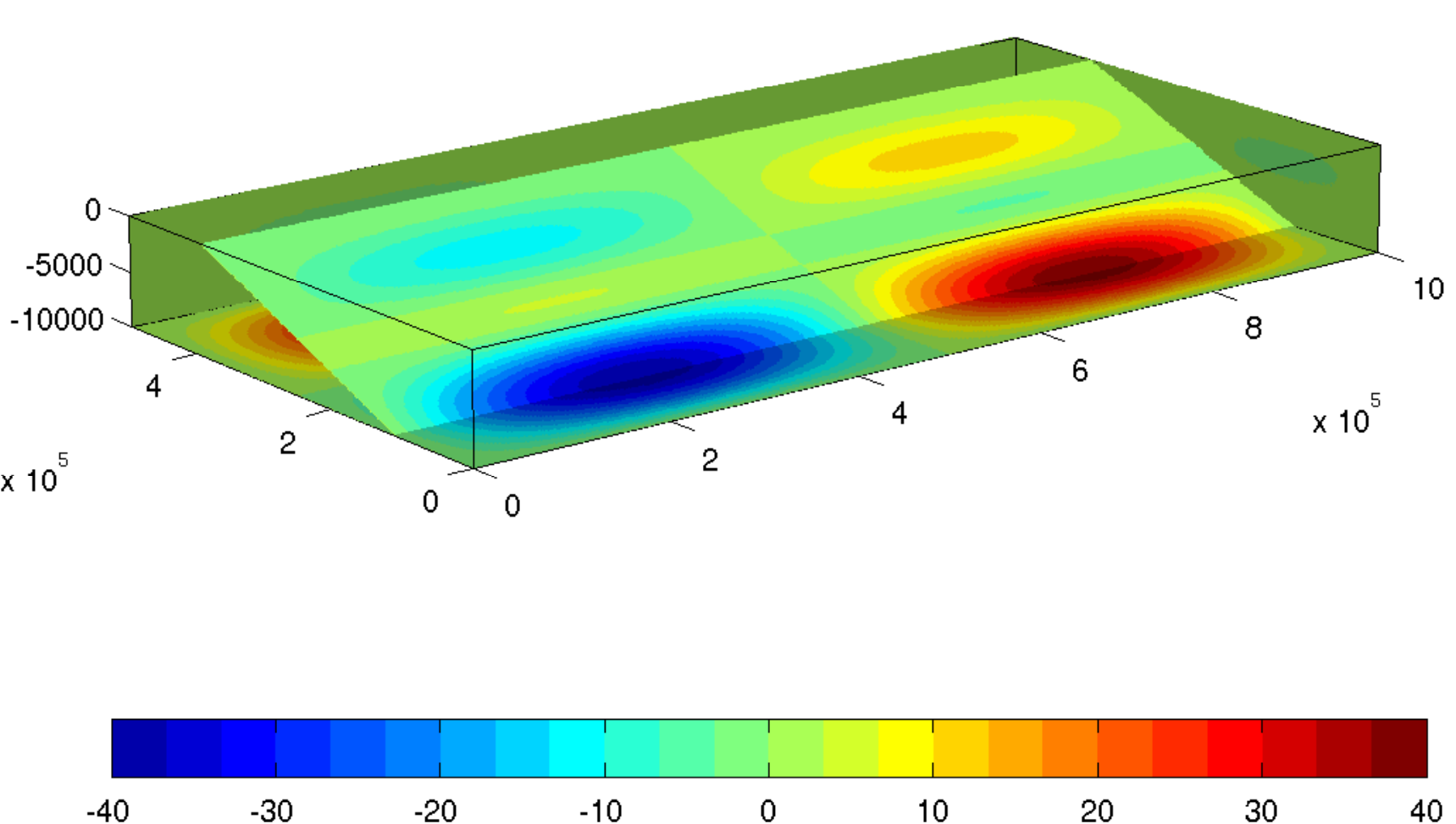}
\caption{The initial state of $\phi$ in the larger domain $\M$}
\label{f1.4}
\end{center}
\end{figure}

\begin{figure}[ht]
   \begin{center}
      \includegraphics[scale=0.8]{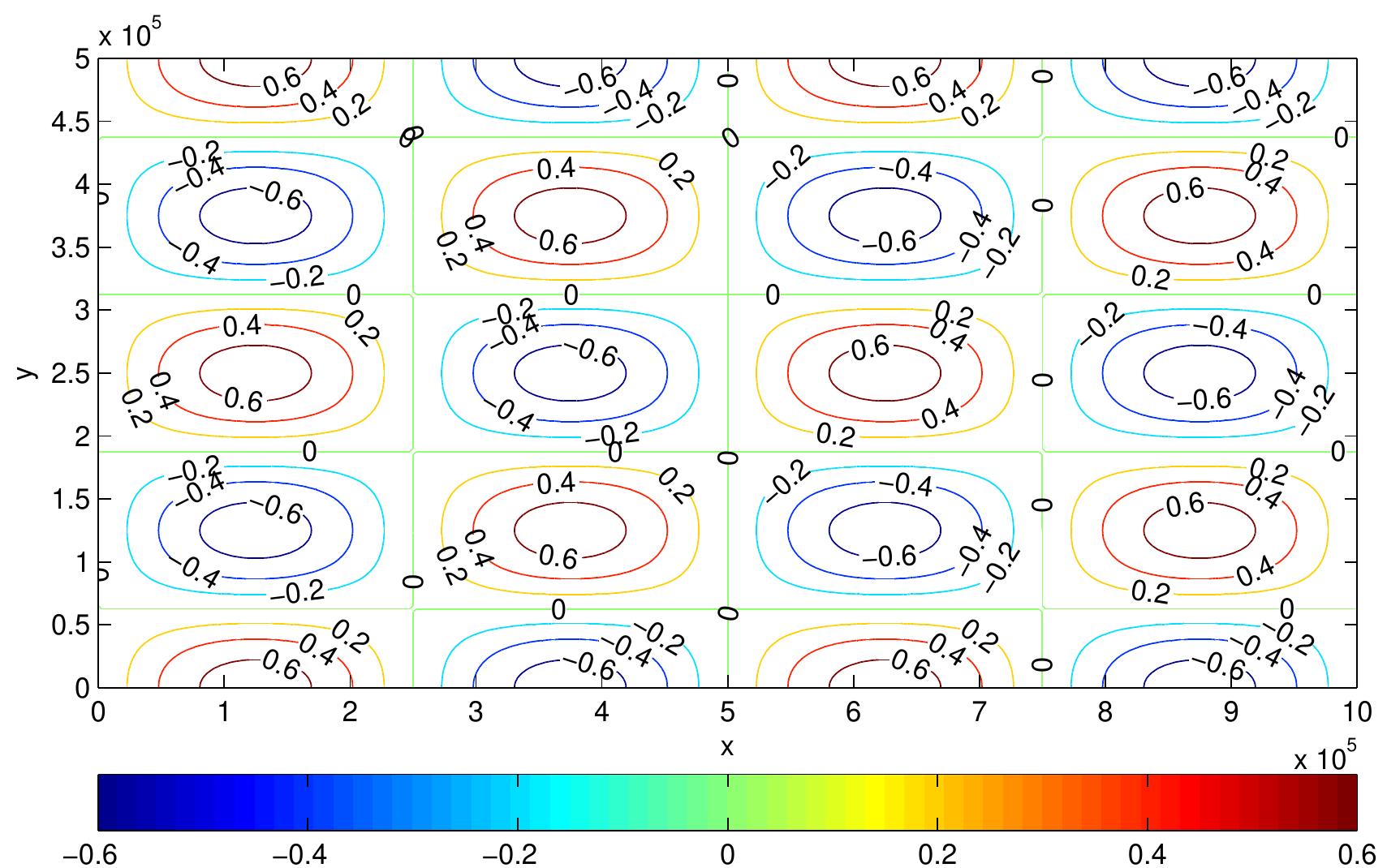}
   \end{center}
   \caption{Contour plot of $u$ at $z=-2500m$, at $t=0$.}
   \label{f1.4a}
\end{figure}

\begin{figure}[ht]
   \begin{center}
      \includegraphics[scale=0.8]{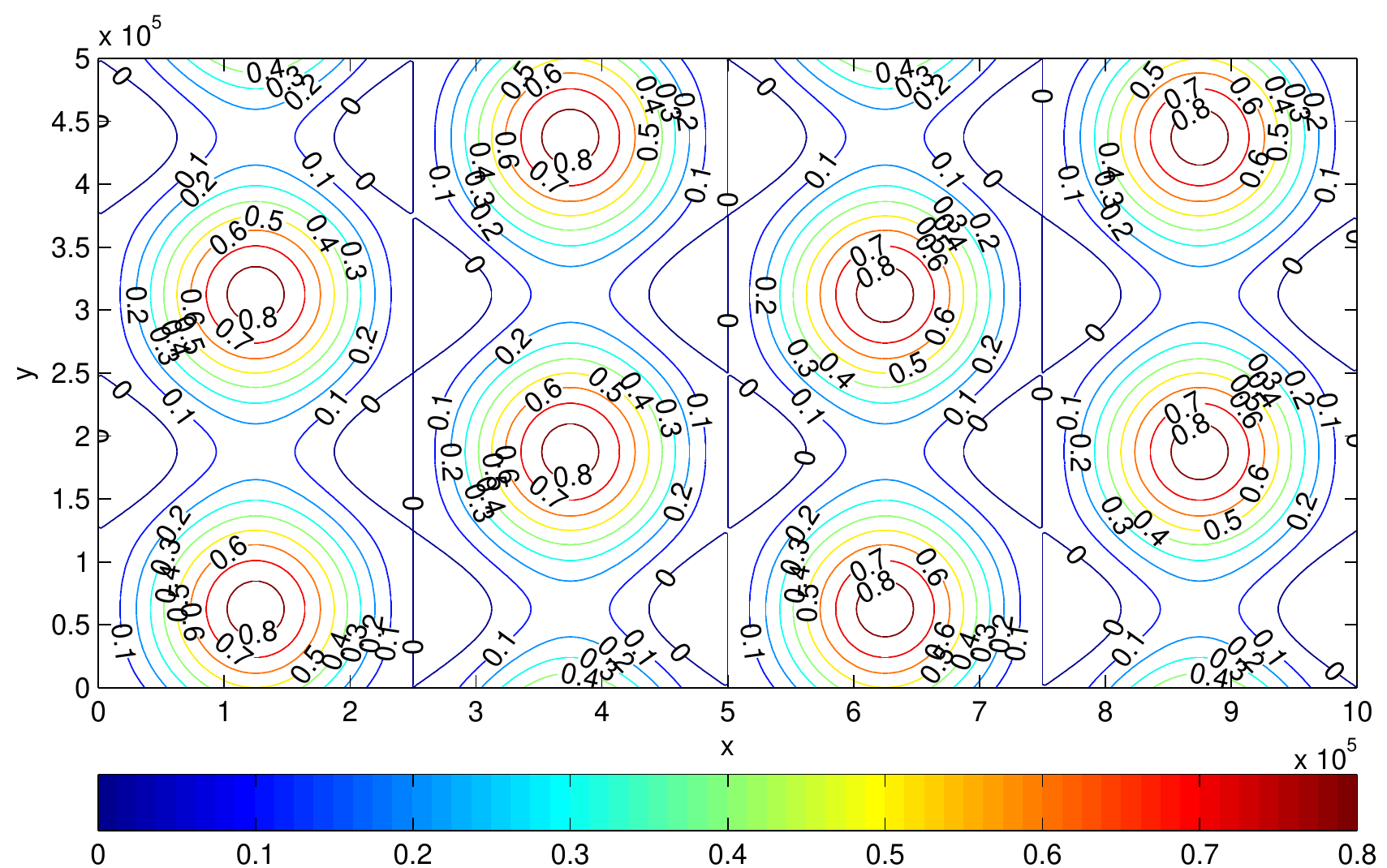}
   \end{center}
   \caption{Contour plot of $v$ at $z=-2500m$, at $t=0$.}
   \label{f1.4b}
\end{figure}

\begin{figure}[ht]
   \begin{center}
      \includegraphics[scale=0.8]{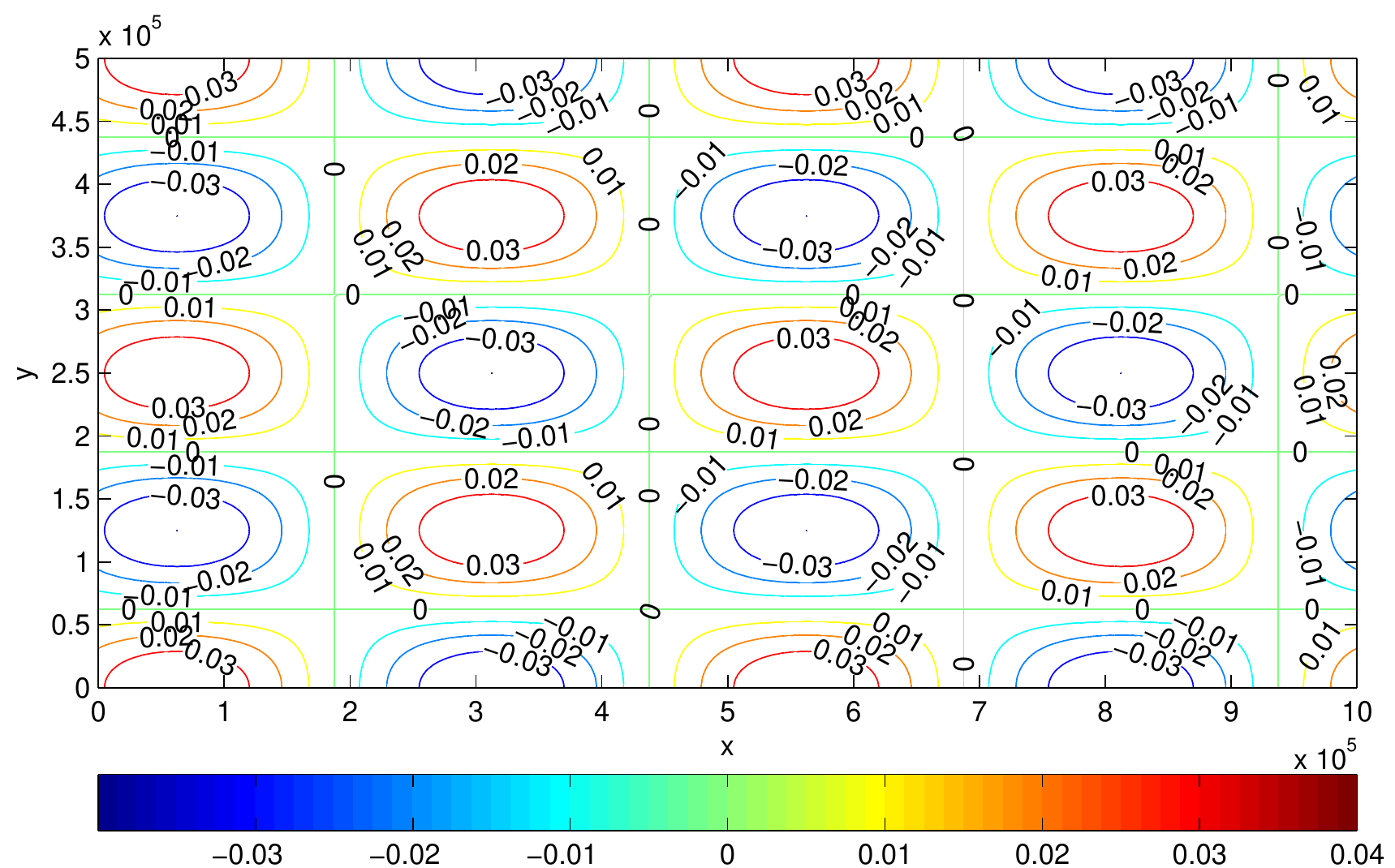}
   \end{center}
   \caption{Contour plot of $w$ at $z=-2500m$, at $t=0$.}
   \label{f1.4c}
\end{figure}

\begin{figure}[ht]
   \begin{center}
      \includegraphics[scale=0.8]{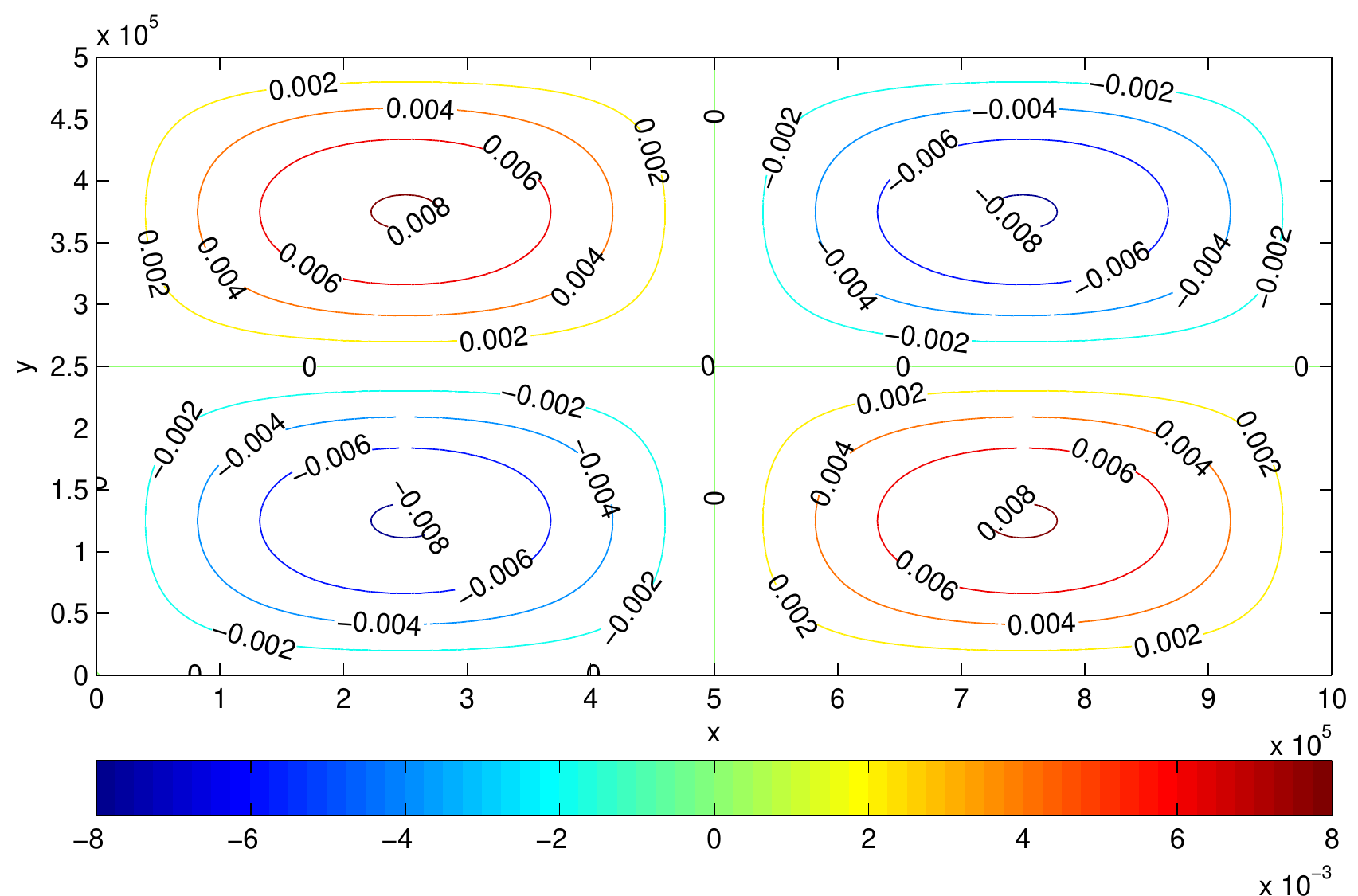}
   \end{center}
   \caption{Contour plot of $\psi$ at $z=-2500m$, at $t=0$.}
   \label{f1.4d}
\end{figure}

\begin{figure}[ht]
   \begin{center}
      \includegraphics[scale=0.8]{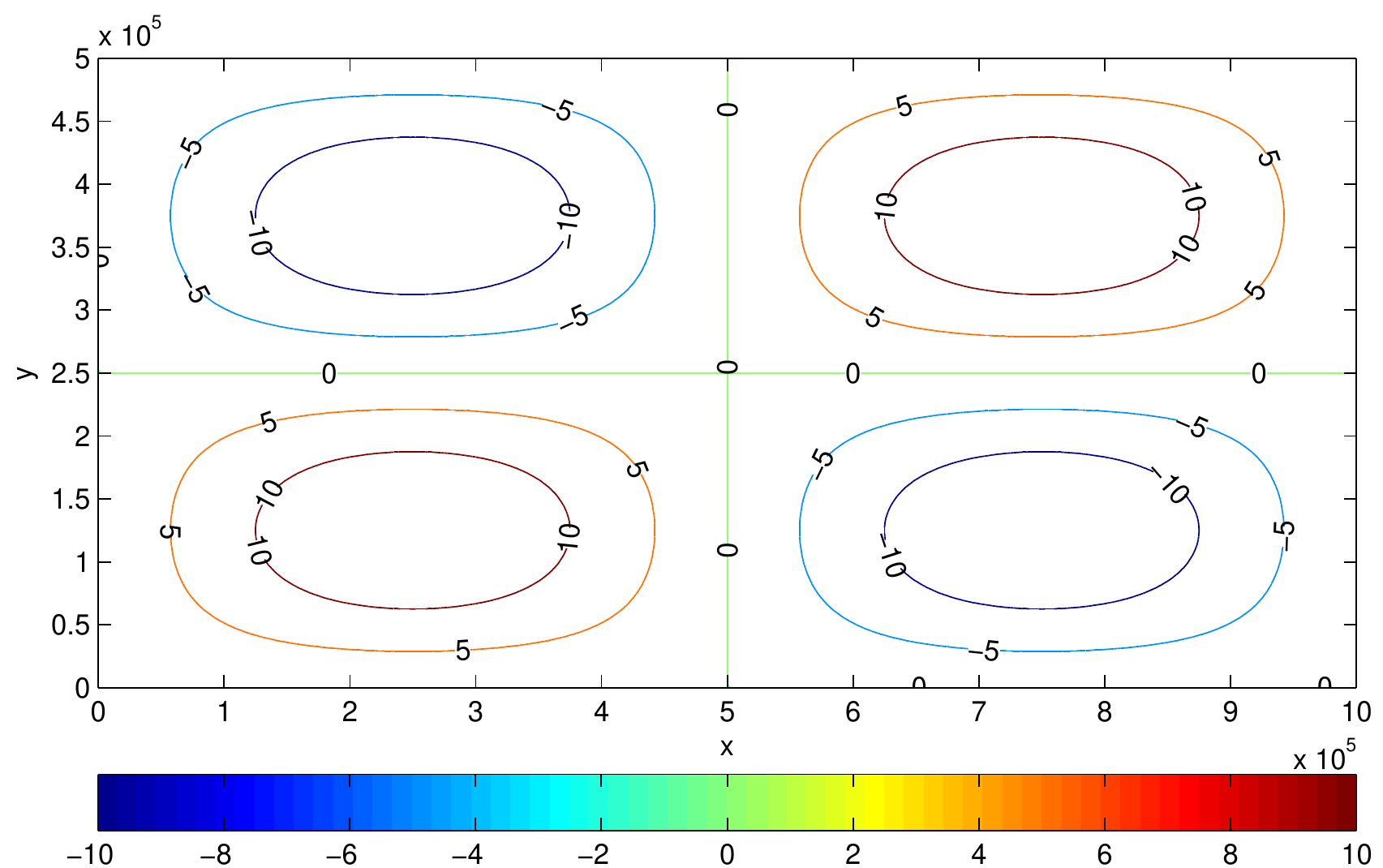}
   \end{center}
   \caption{Contour plot of $\phi$ at $z=-2500m$, at $t=0$.}
   \label{f1.4e}
\end{figure}

\begin{figure}
\begin{center}
\includegraphics[scale=0.8]{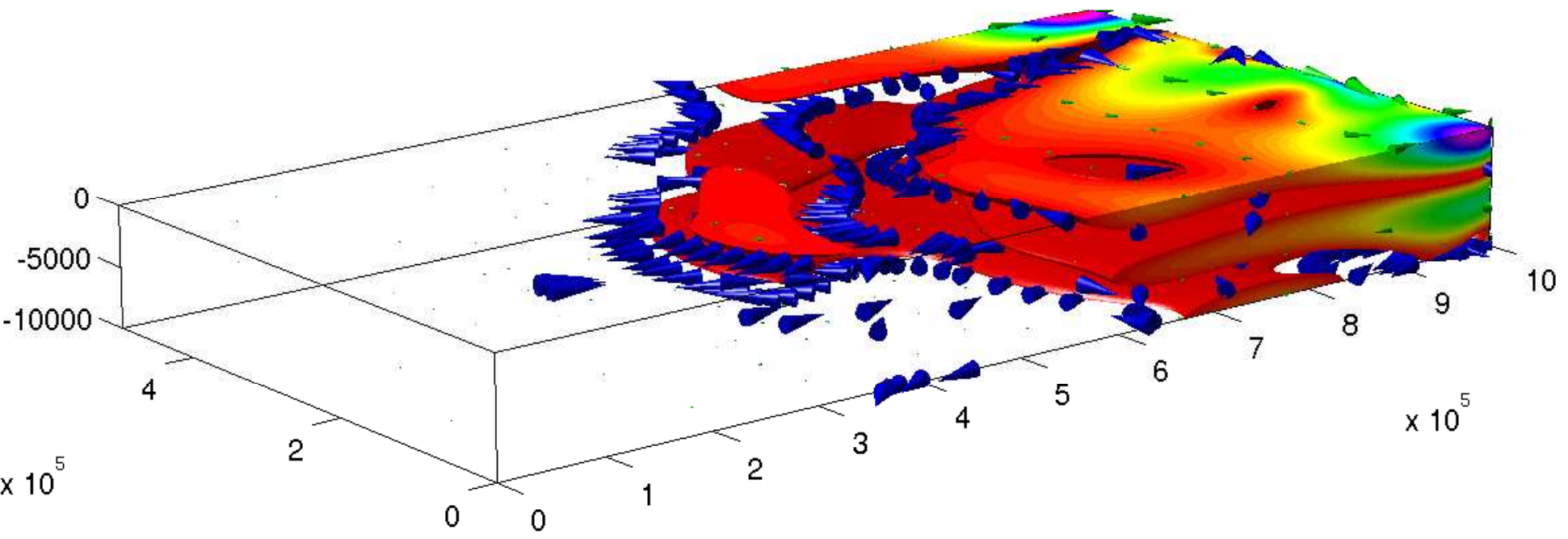}
\end{center}
\caption{The velocity field with cone plot in the larger domain $\M$ at $t = T$}
\label{f1.5}
\end{figure}

\begin{figure}
\begin{center}
\includegraphics[scale=0.9]{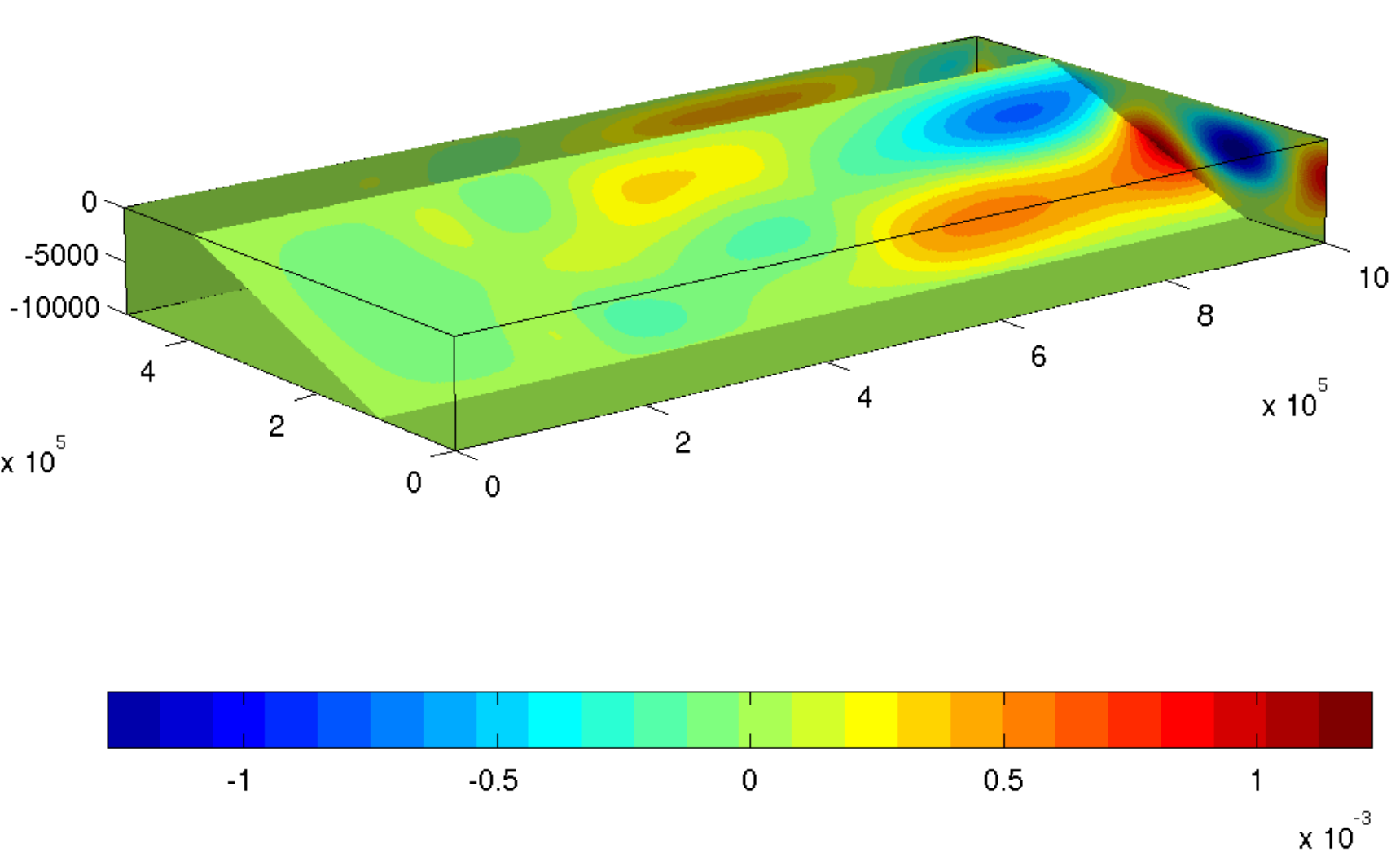}
\end{center}
\caption{The state of $\psi$ in the larger domain $\M$ at $t = T$. }
\label{f1.6}
\end{figure}

\begin{figure}
\begin{center}
\includegraphics[scale=0.9]{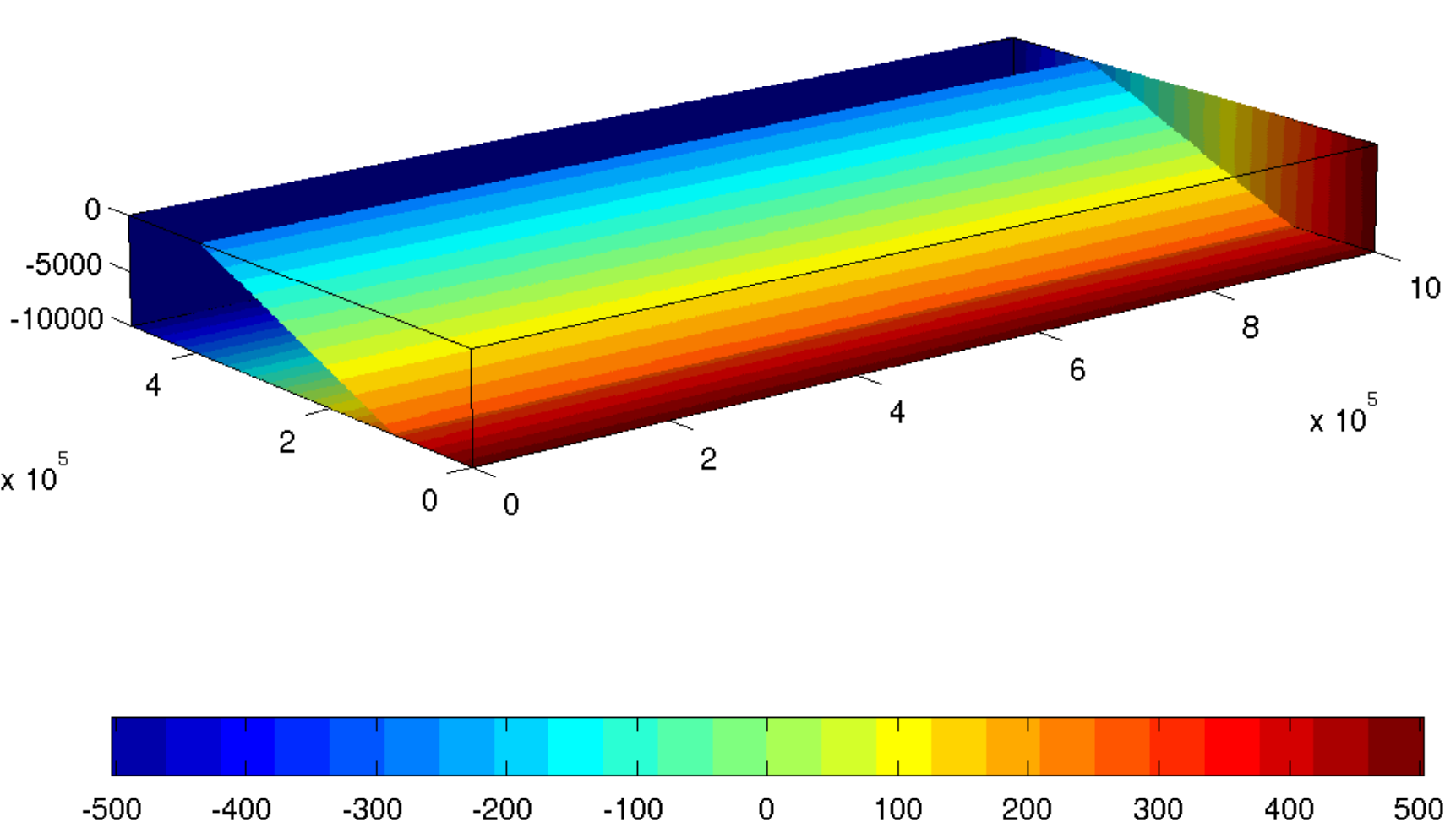}
\end{center}
\caption{The state of $\phi$ in the larger domain $\M$ at $t = T$.}
\label{f1.7}
\end{figure}

\begin{figure}[ht]
   \begin{center}
      \includegraphics[scale=0.8]{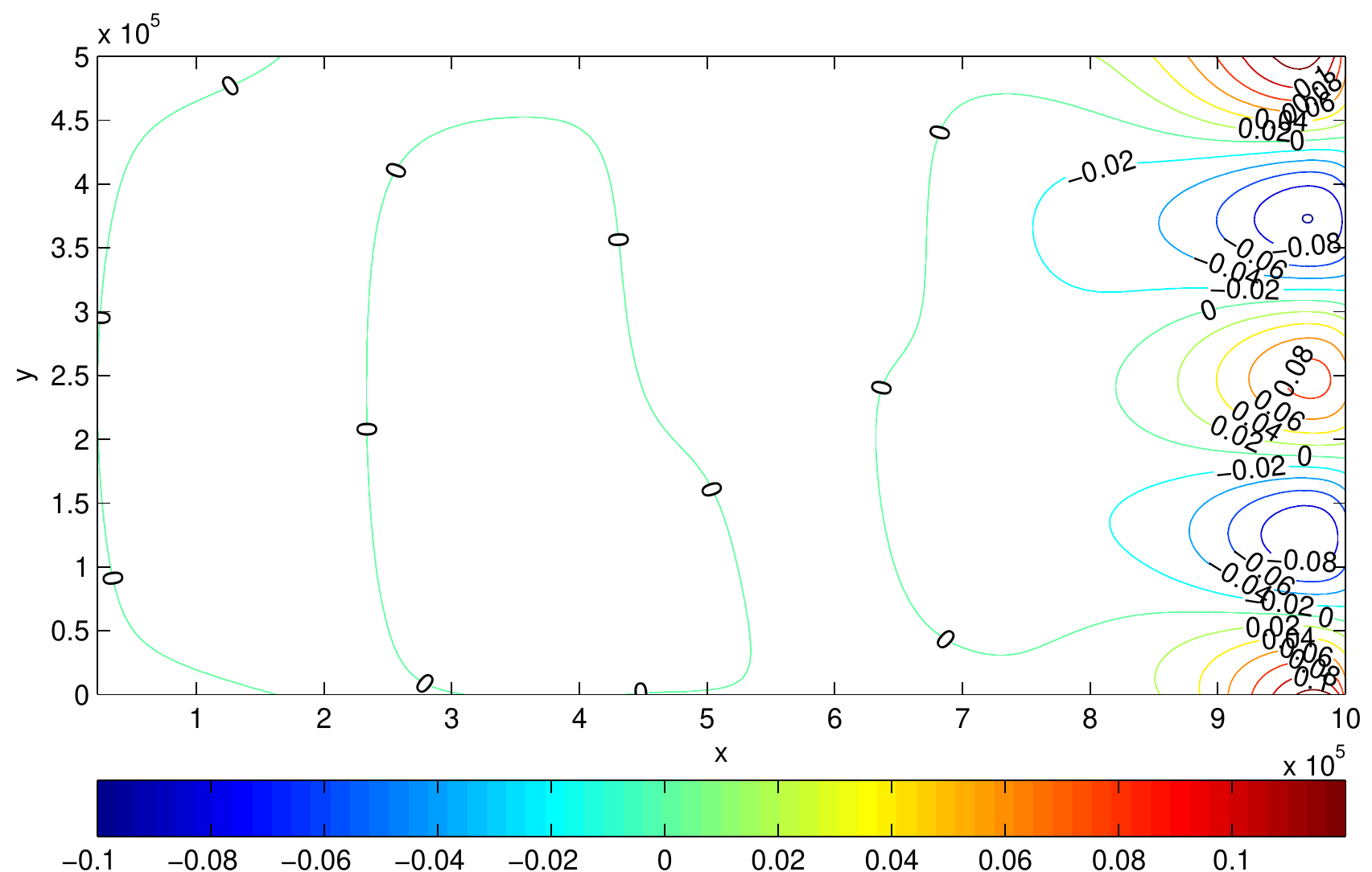}
   \end{center}
   \caption{Contour plot of $u$ at $z=-2500m$, at $t=T$.}
   \label{f1.7a}
\end{figure}

\begin{figure}[ht]
   \begin{center}
      \includegraphics[scale=0.8]{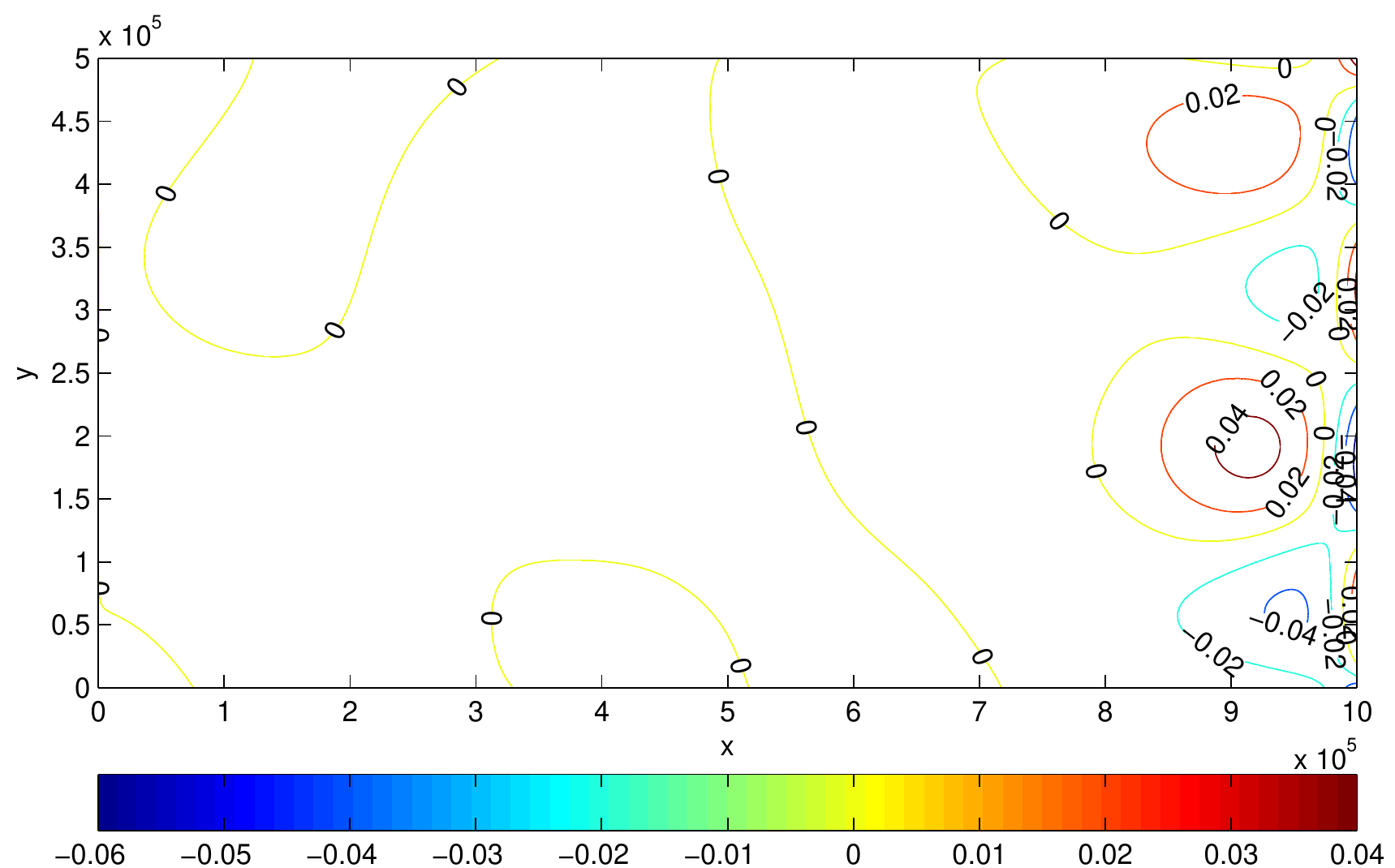}
   \end{center}
   \caption{Contour plot of $v$ at $z=-2500m$, at $t=T$.}
   \label{f1.7b}
\end{figure}

\begin{figure}[ht]
   \begin{center}
      \includegraphics[scale=0.8]{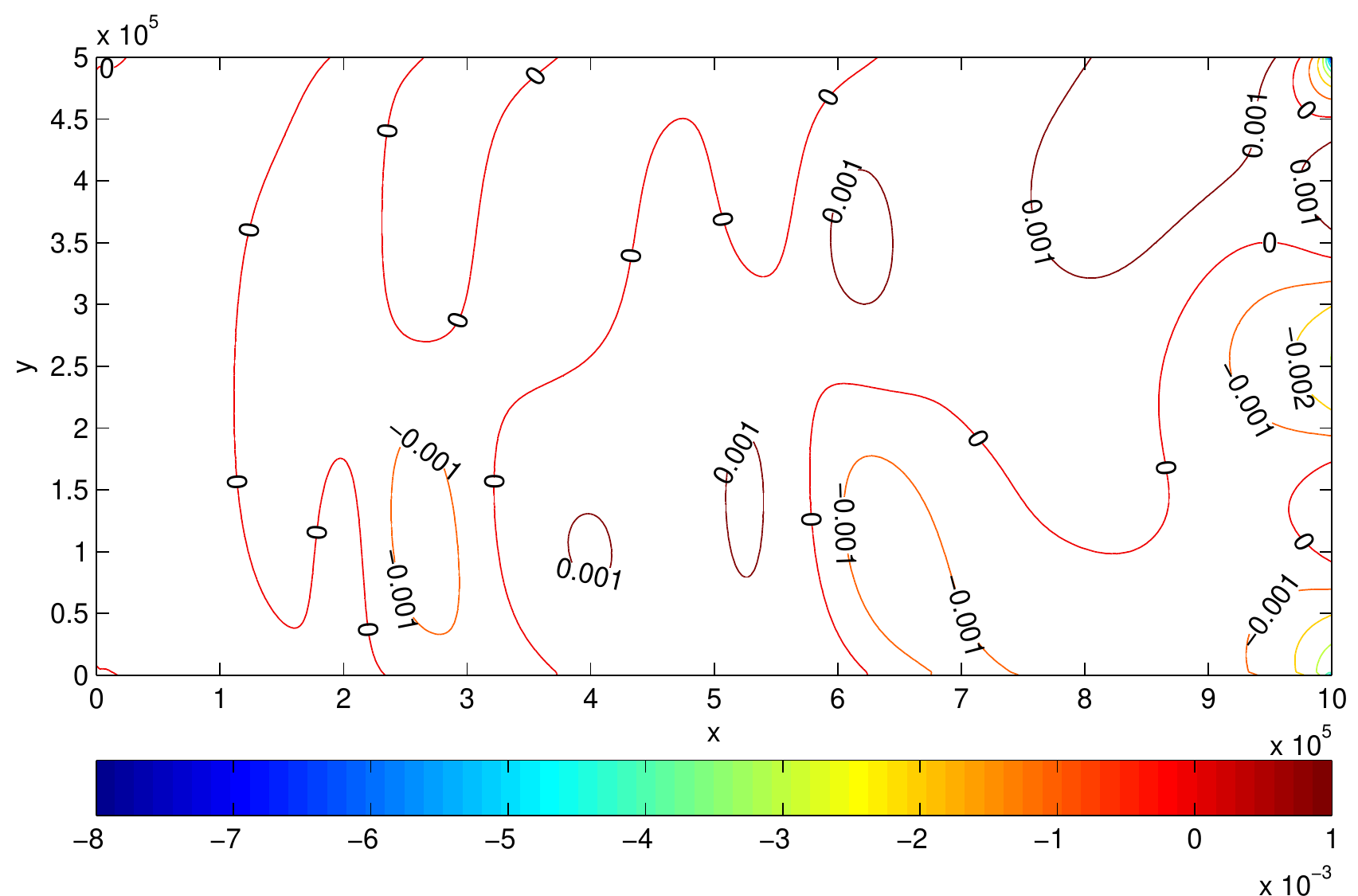}
   \end{center}
   \caption{Contour plot of $w$ at $z=-2500m$, at $t=T$.}
   \label{f1.7c}
\end{figure}

\begin{figure}[ht]
   \begin{center}
      \includegraphics[scale=0.8]{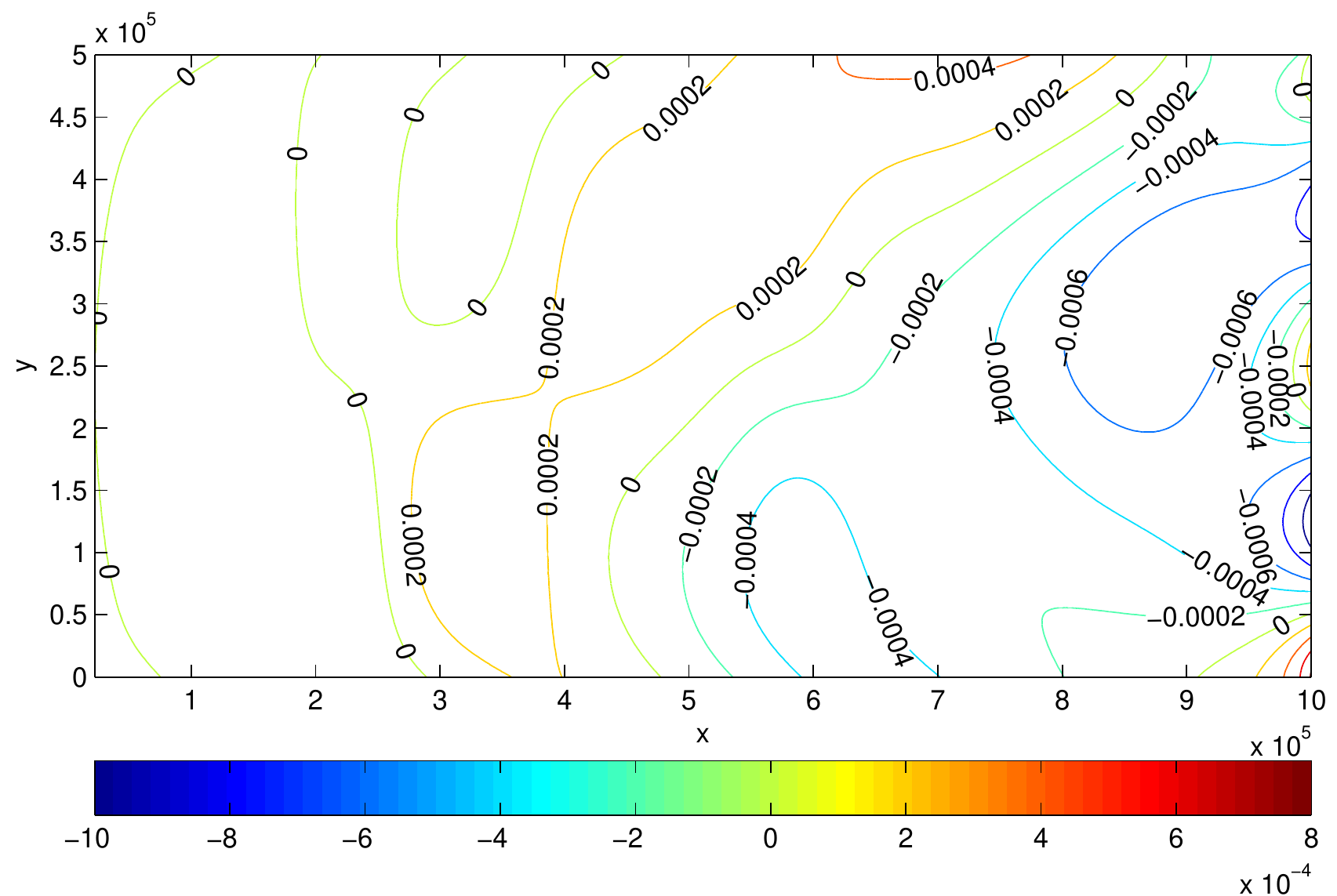}
   \end{center}
   \caption{Contour plot of $\psi$ at $z=-2500m$, at $t=T$.}
   \label{f1.7d}
\end{figure}

\begin{figure}[ht]
   \begin{center}
      \includegraphics[scale=0.8]{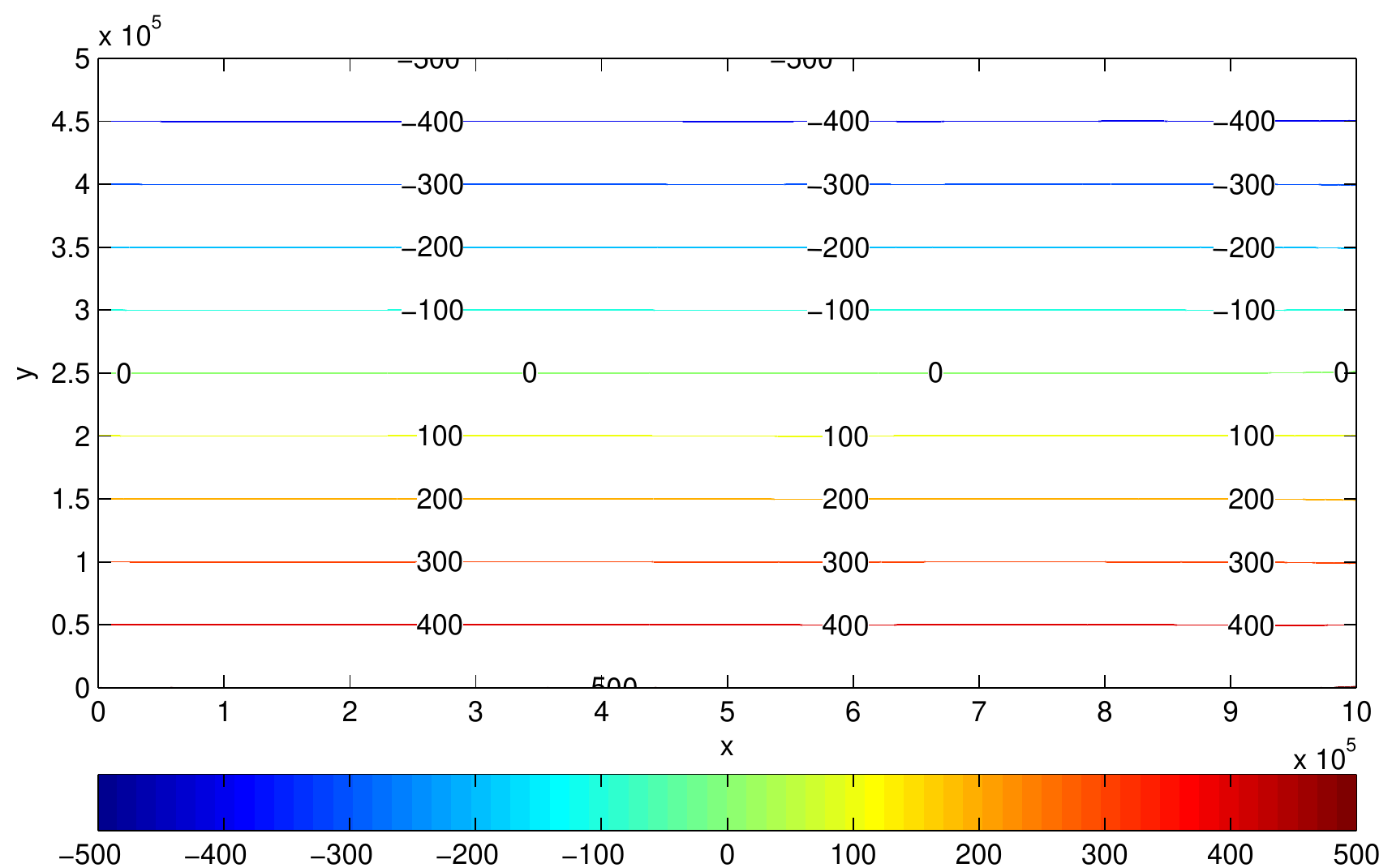}
   \end{center}
   \caption{Contour plot of $\phi$ at $z=-2500m$, at $t=T$.}
   \label{f1.7e}
\end{figure}

The simulation results over the larger domain $\M$ are plotted
in Figures \ref{f1.2} to \ref{f1.7e}.
Figure \ref{f1.2} is the cone plot with isosurface of the initial state of the velocity field, and Figures \ref{f1.3} and \ref{f1.4} are the slice--plane plots of the initial state of $\phi$ and $\psi$.
Figures \ref{f1.4a} to \ref{f1.4e} are the contour plots of $u$, $v$,
$w$, $\psi$ and $\phi$, respectively, on the plane $z=-2,500m$, at
$t=0$.

Figure \ref{f1.5} is the cone plot with isosurface of the velocity field
at the final time $t =T$, and Figures \ref{f1.6} and \ref{f1.7} are
the slice--plane plots of the state of $\phi$ and $\psi$
at the final time $t = T$.
Figures \ref{f1.7a} to \ref{f1.7e} are the contour plots of $u$, $v$,
$w$, $\psi$ and $\phi$, respectively, on the plane $z=-2,500m$, at
$t=T$.
\clearpage

\subsection{Simulation in the middle-half domain}\label{s4.2}
As explained in the Introduction, we next do simulations on the middle half domain $\M_1$ of $\M$, $\M_1 = \M_1'\times (-H,0)$, $\M_1' = (\frac{L_1}{4},\frac{3\, L_1}{4})\times (\frac{L_2}{4},\frac{3\, L_2}{4}).$
The \emph{boundary values of the unknown functions
$u,\, v,\, w,\, \phi,$ and $\psi$ are inferred from the previous simulation}. More specifically, the boundary conditions are, for the zeroth mode ($n=0$),
\begin{equation}\label{e7.5}
\begin{cases}
u_0(L_1/4, y_j, t_k) = u^l_0(L_1/4,y_j,t_k), \\
u_0(3L_1/4,y_j,t_k) = u^l_0(3L_1/4,y_j,t_k),\\
v_0(L_1/4,y_j,t_k) = v^l_0(L_1/4,y_j,t_k),\\
v_0(x_i,L_2/4,t_k) = v^l_0(x_i,L_2/4,t_k),\\
v_0(x_i,3L_2/4,t_k) = v^l_0(x_i,3L_2/4,t_k).
\end{cases}
\end{equation}
For the subcritical modes ($ 1\le n \le n_c$),
\begin{equation}\label{e7.6}
\begin{cases}
\xi_n(L_1/4,y_j,t_k) = \xi^l_n(L_1/4,y_j,t_k),\\
v_n(L_1/4,y_j,t_k) = v^l_n(L_1/4,y_j,t_k),\\
\eta_n(3L_1/4,y_j,t_k) = \eta^l_n(3L_1/4,y_j,t_k),\\
\alpha_n(x_i,3L_2/4,t_k) = \alpha^l_n(x_i,3L_2/4,t_k), \\
\beta_n(x_i,L_2/4,t_k) = \beta^l_n(x_i, L_2/4, t_k),
\end{cases}
\end{equation}
and for the supercritical modes ($ n > n_c$),
\begin{equation}\label{e7.7}
\begin{cases}
\xi_n(L_1/4,y_j,t_k) = \xi^l_n(L_1/4,y_j,t_k),\\
v_n(L_1/4,y_j,t_k) = v^l_n(L_1/4,y_j,t_k),\\
\eta_n(L_1/4,y_j,t_k) = \eta^l_n(L_1/4,y_j,t_k),\\
\alpha_n(x_i,3L_2/4,t_k) = \alpha^l_n(x_i,3L_2/4,t_k), \\
\beta_n(x_i,L_2/4,t_k) = \beta^l_n(x_i, L_2/4, t_k).
\end{cases}
\end{equation}
In the above, $x_i$, $y_j$ and $t_k$ denote the discrete grid points in space and time. The superscript $l$ denotes the previous simulation in the larger domain $\M$. In this simulation, we take 200 segments in the $x$-direction, and 100 segments in the $y$ direction.

\begin{figure}
\begin{center}
\includegraphics[scale=0.8]{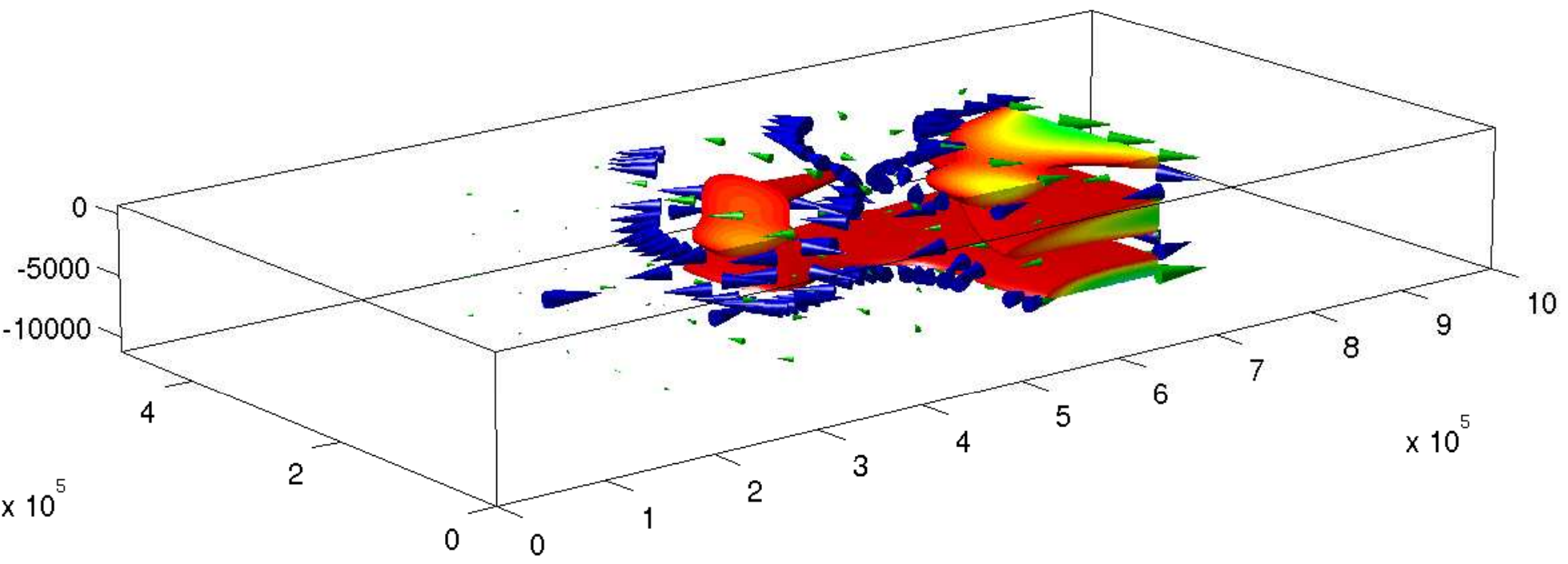}
\end{center}
\caption{The velocity field with cone plot in the middle half domain $\M_1$ at $t = T$}
\label{f1.8}
\end{figure}

\begin{figure}
\begin{center}
\includegraphics[scale=0.9]{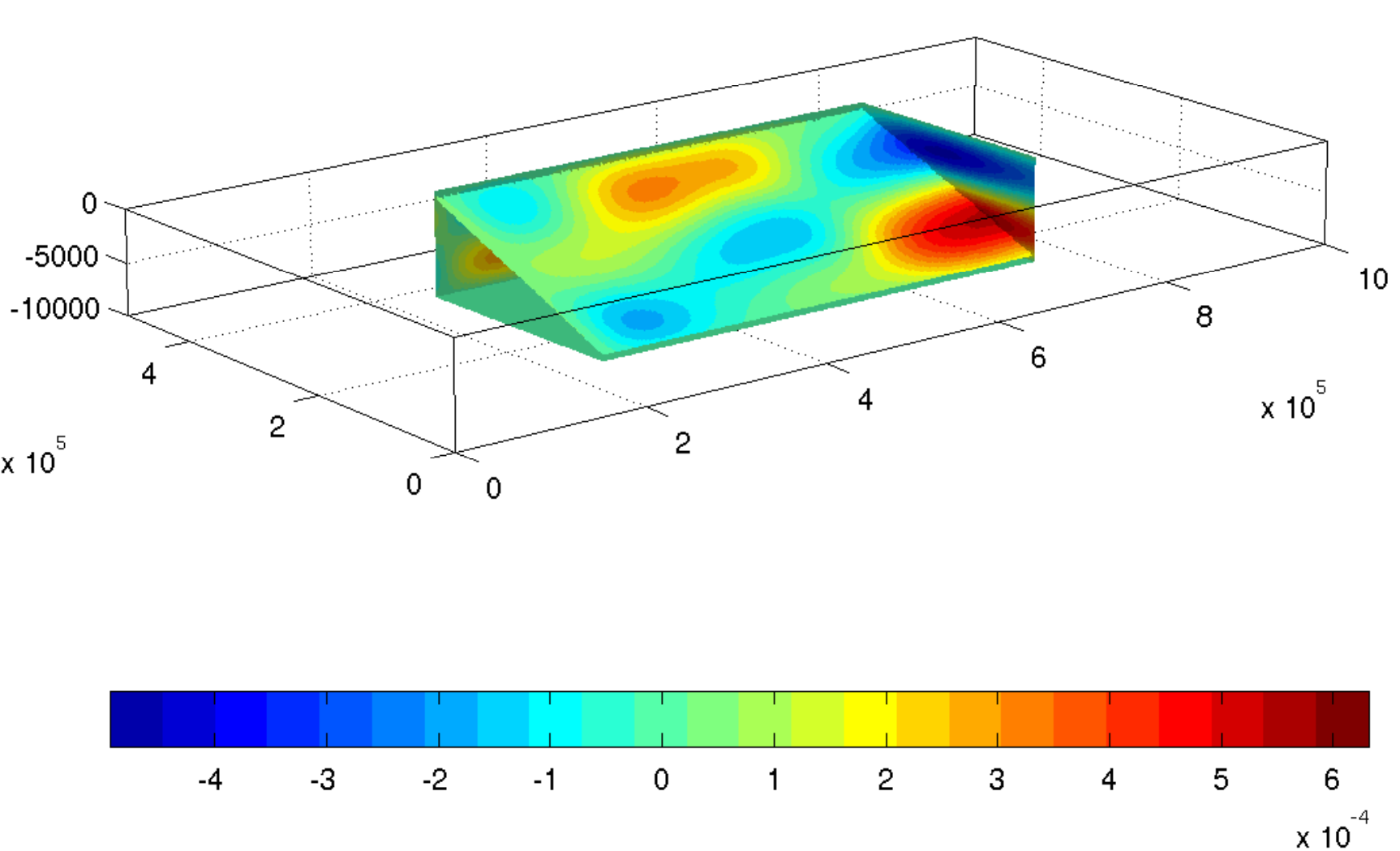}
\end{center}
\caption{The state of $\psi$ in the middle half domain $\M_1$ at $t = T$. }
\label{f1.9}
\end{figure}

\begin{figure}
\begin{center}
\includegraphics[scale=0.9]{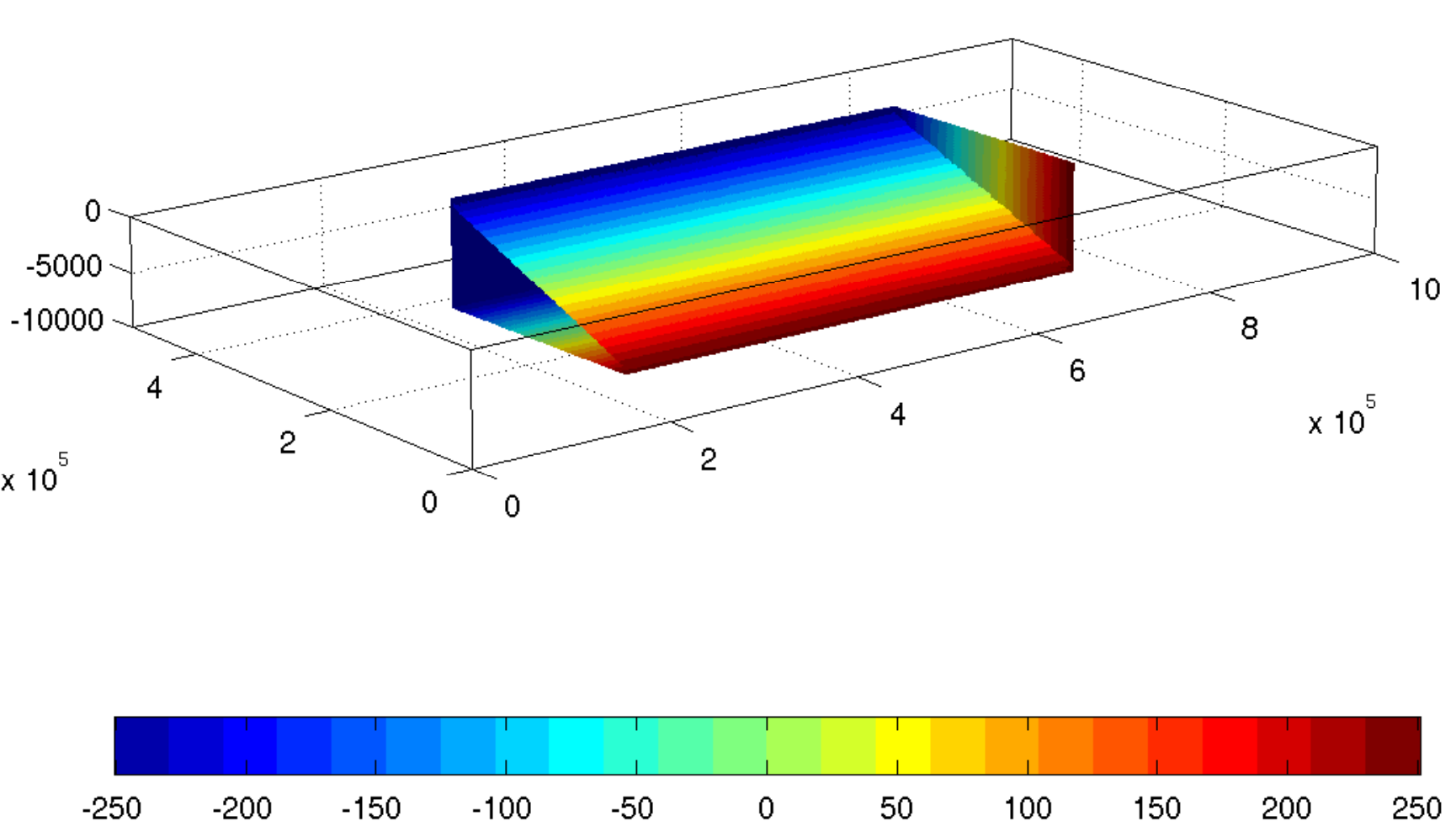}
\end{center}
\caption{The state of $\phi$ in the middle half domain $\M_1$ at $t = T$.}
\label{f1.10}
\end{figure}

\begin{figure}[ht]
   \begin{center}
      \includegraphics[scale=0.8]{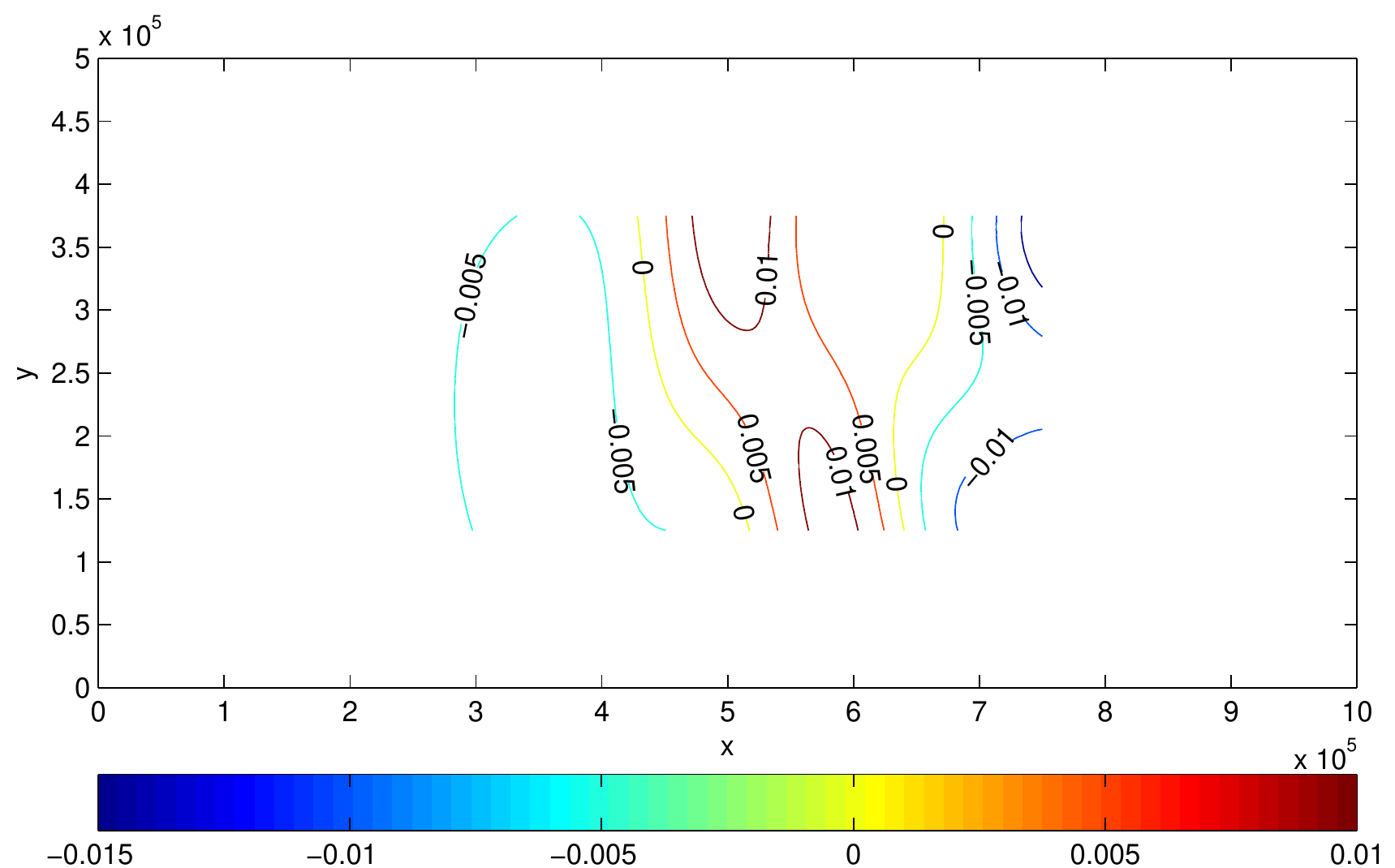}
   \end{center}
   \caption{Contour plot of $u$ at $z=-2500m$, at $t=T$.}
   \label{f1.10a}
\end{figure}

\begin{figure}[ht]
   \begin{center}
      \includegraphics[scale=0.8]{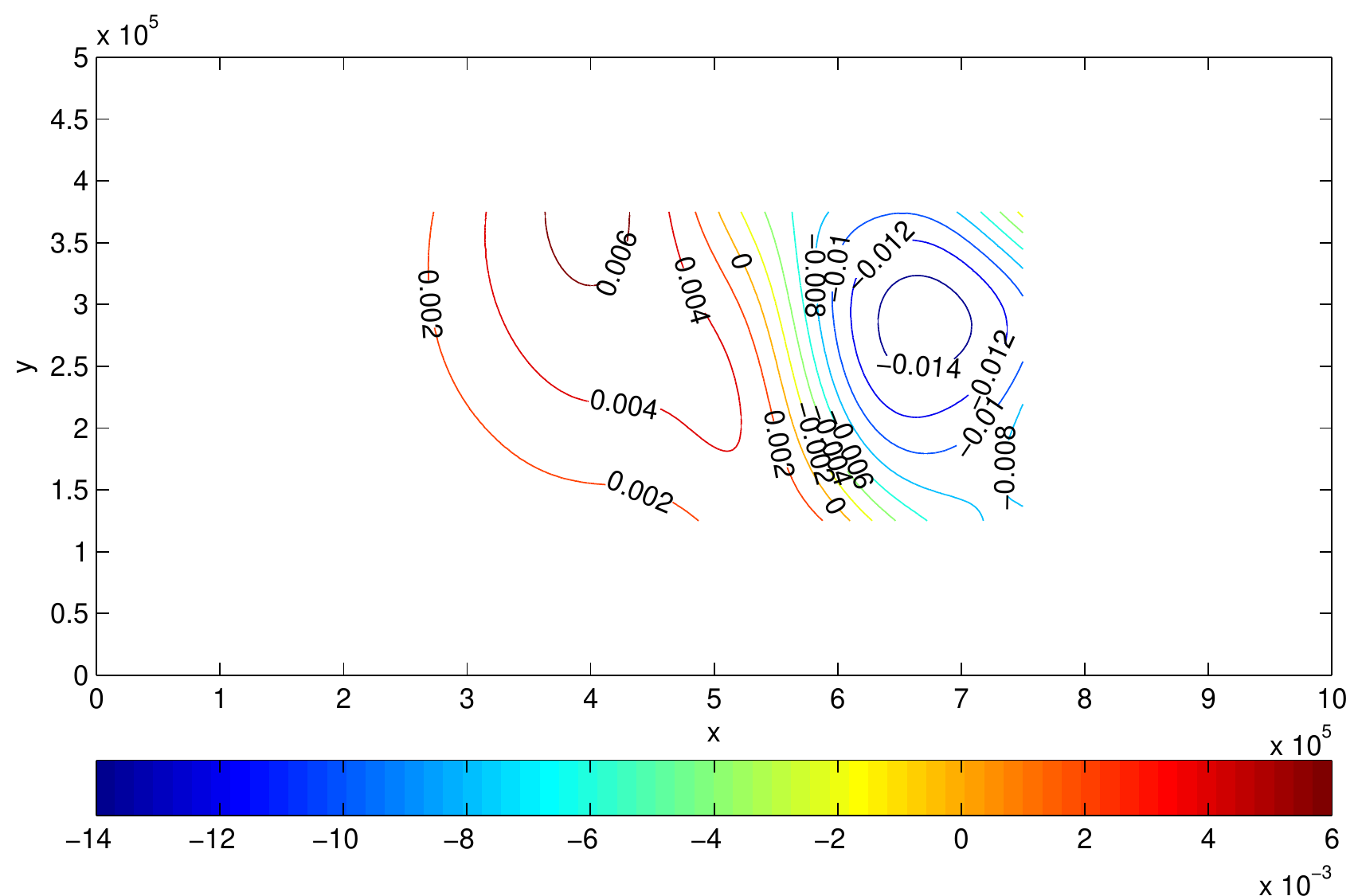}
   \end{center}
   \caption{Contour plot of $v$ at $z=-2500m$, at $t=T$.}
   \label{f1.10b}
\end{figure}

\begin{figure}[ht]
   \begin{center}
      \includegraphics[scale=0.8]{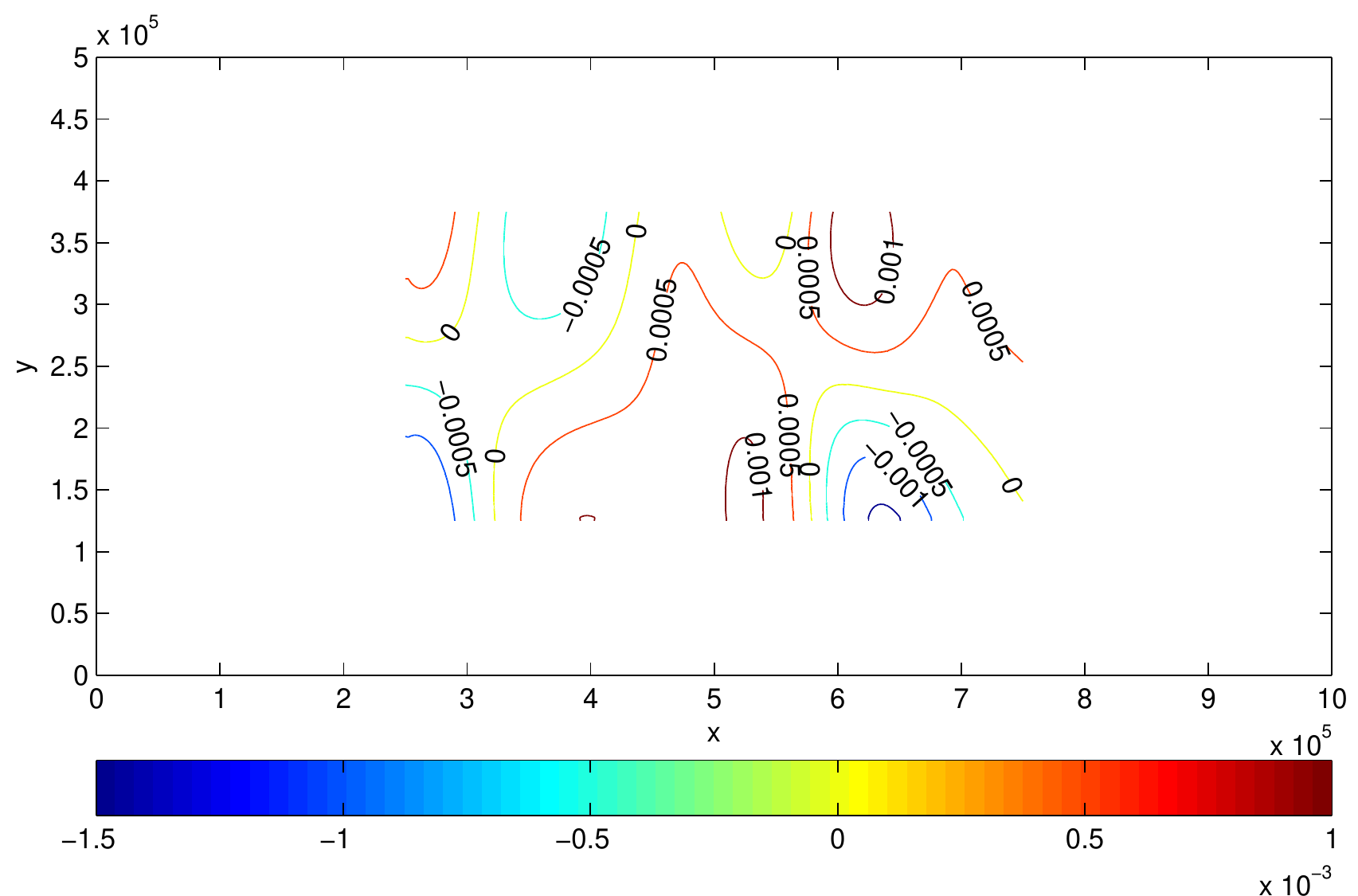}
   \end{center}
   \caption{Contour plot of $w$ at $z=-2500m$, at $t=T$.}
   \label{f1.10c}
\end{figure}

\begin{figure}[ht]
   \begin{center}
      \includegraphics[scale=0.8]{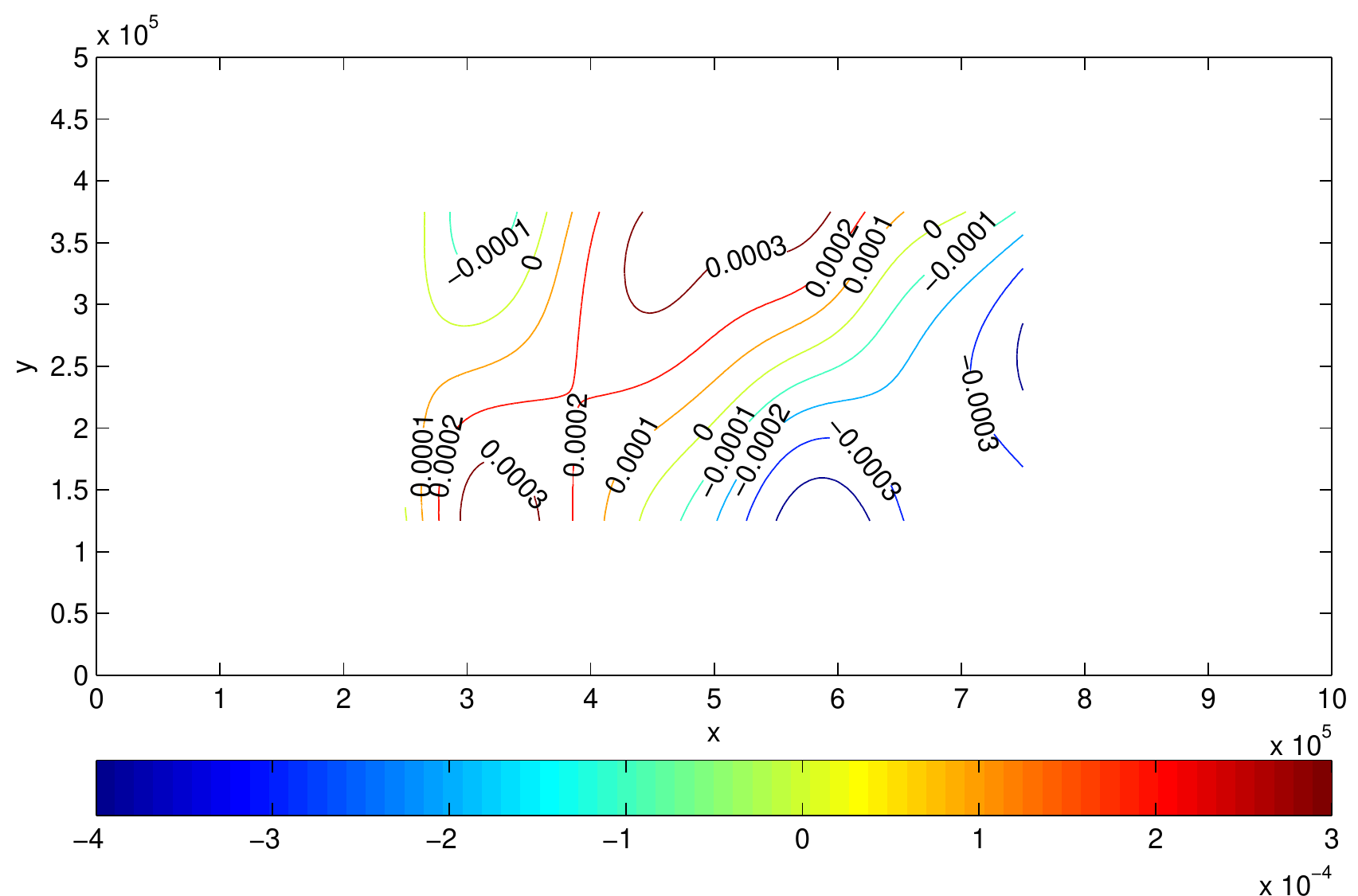}
   \end{center}
   \caption{Contour plot of $\psi$ at $z=-2500m$, at $t=T$.}
   \label{f1.10d}
\end{figure}

\begin{figure}[ht]
   \begin{center}
      \includegraphics[scale=0.8]{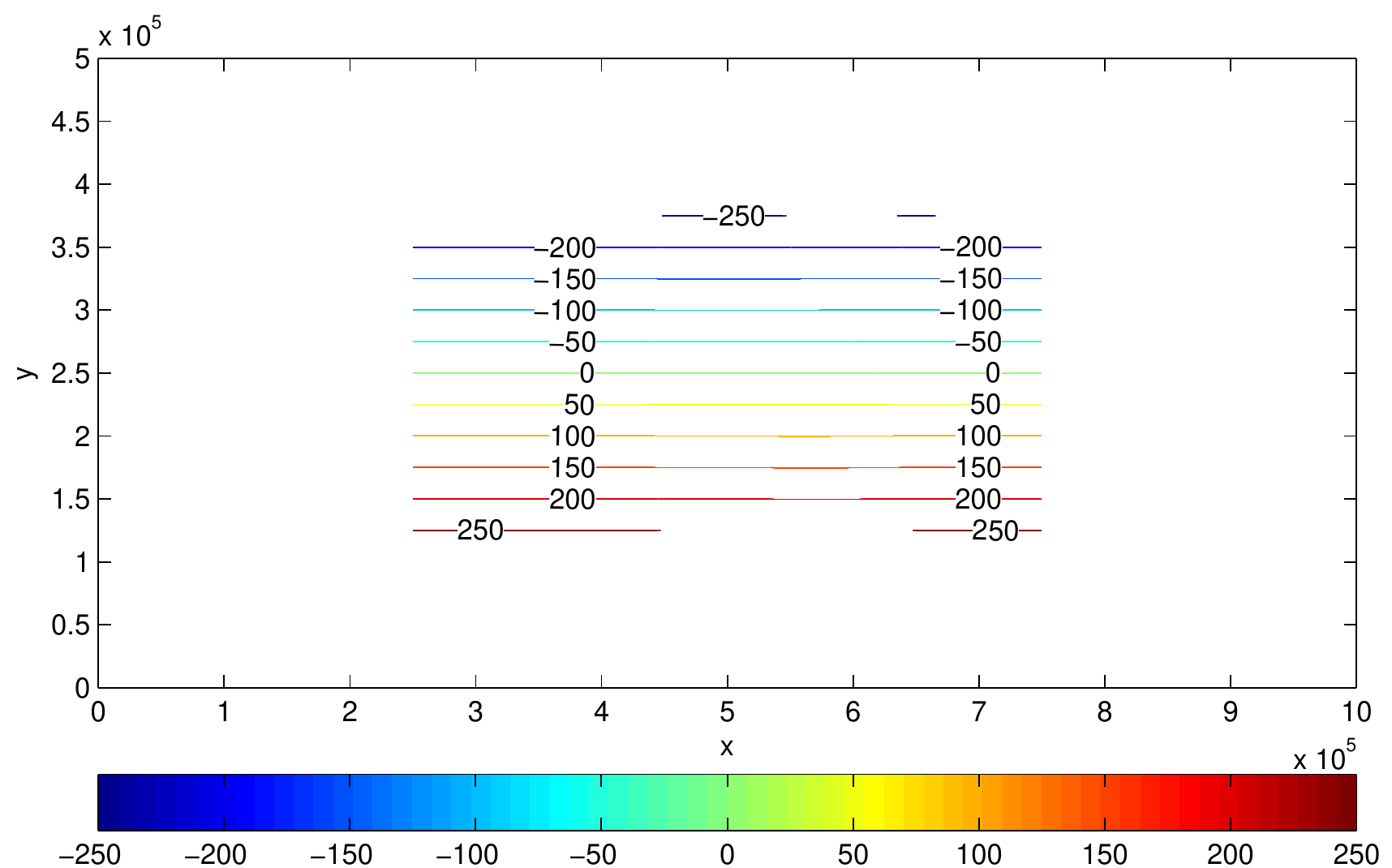}
   \end{center}
   \caption{Contour plot of $\phi$ at $z=-2500m$, at $t=T$.}
   \label{f1.10e}
\end{figure}

The simulation results over the middle half domain $\M_1$
are plotted in Figures \ref{f1.8} to \ref{f1.10e}.
Figure \ref{f1.8} is the cone plot with isosurface of
the velocity field
in the middle half domain $\M_1$ at the final time $t =T$, and
Figures \ref{f1.9} and \ref{f1.10} are the slice--plane plots of
the state of $\phi$ and $\psi$
in the middle half domain $\M_1$ at the final time $t = T$.
Figures \ref{f1.10a} and \ref{f1.10e} are the contour plots of $u$, $v$,
$w$, $\psi$ and $\phi$, respectively, on the plane $z=-2,500m$ restricted
to the middle half domain $\M_1$, at
$t=T$.


\subsection{Comparison}\label{s4.3}
In this subsection, \emph{we compare these two distinct simulations namely the results of the simulations on the
larger domain $\M$ restricted to the middle half domain $\M_1$
and the results of the simulations on the middle half domain $\M_1$, as obtained
by the second simulation above}.

Let $u^{ext}, \, v^{ext},\, w^{ext},\,\phi^{ext}$ and $\psi^{ext}$ be the numerical approximations of the variables $u,\, v,\, w,\, \phi$, and $\psi$ on the larger domain $\M$, respectively, and  $u^{int}, \, v^{int},\, w^{int},\,\phi^{int}$ and $\psi^{int}$ be the numerical approximations of the variables $u,\, v,\, w,\, \phi$, and $\psi$ on the middle half domain $\M_1$, respectively.
In Figures \ref{f1.11}-\ref{f1.15}, we plot the evolution of the unknowns,
and of their relative errors (see below), in both the $L^2$ and $L^\infty$ norms.
The relative errors are defined as
$\dfrac{||u^{int}-u^{ext}|_{\M_1}||_{L^p}}{||u^{ext}||_{L^p}}$,
etc.~where $p = 2, \infty$.
\begin{figure}
   \begin{center}
      \includegraphics[scale=0.8]{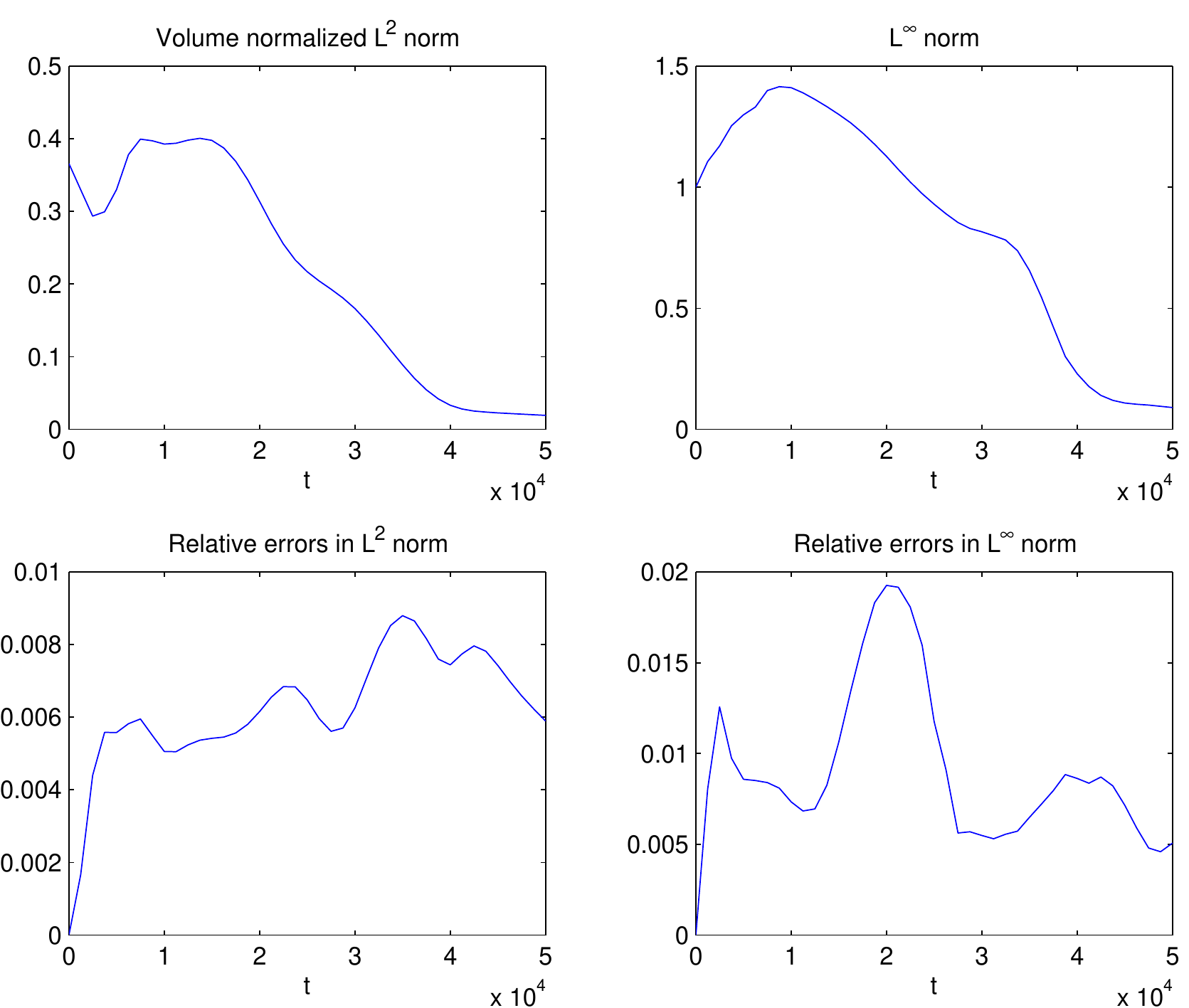}
      \caption{Top row: evolution of  the solution $u$ in $L^2$ and $L^\infty$ norms.
      Bottom row: evolution of the relative errors for $u$ in $L^2$ and $L^\infty$ norms.}
      \label{f1.11}
   \end{center}
\end{figure}

\begin{figure}
   \begin{center}
      \includegraphics[scale=0.8]{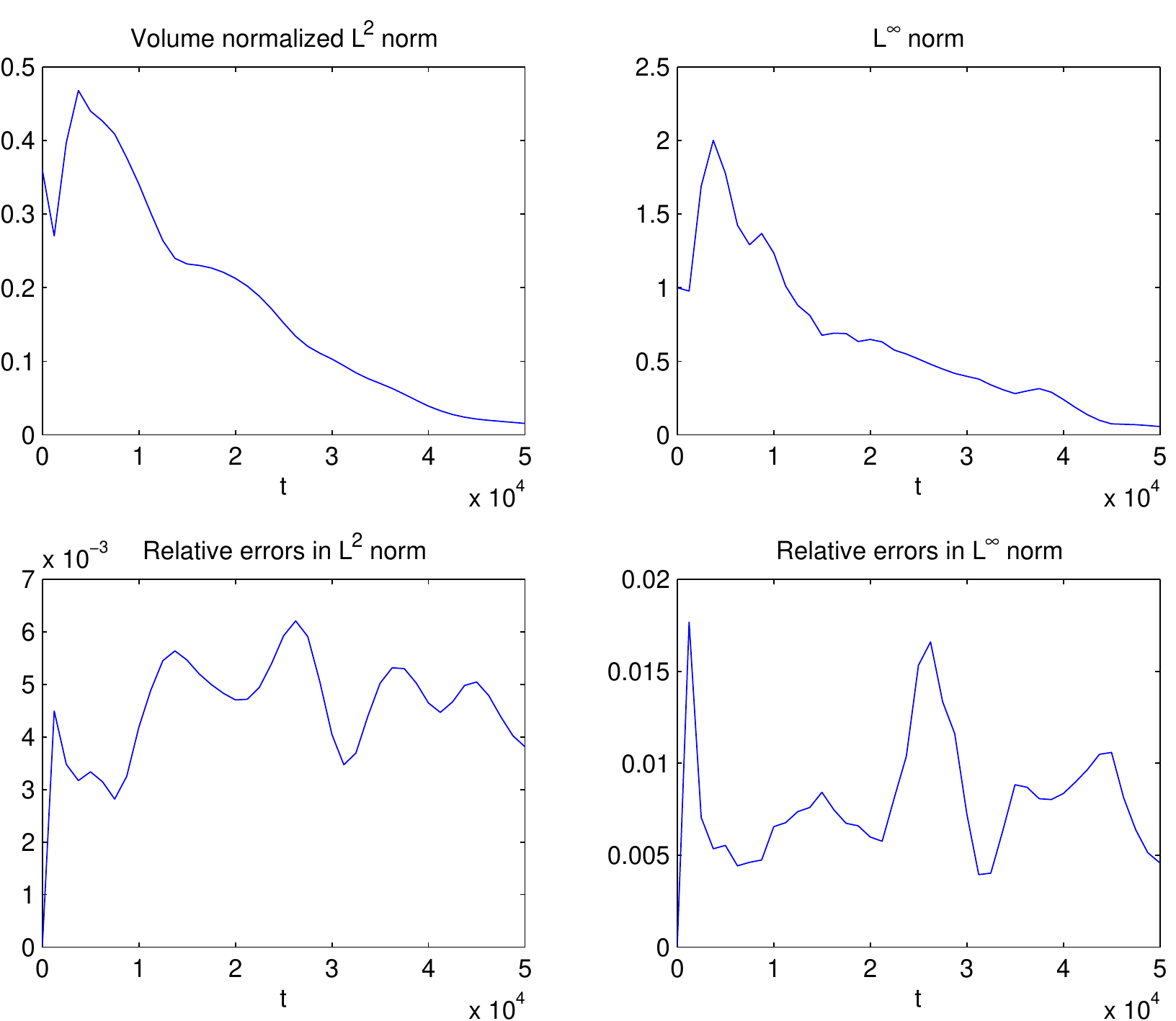}
      \caption{Top row: evolution of  the solution $v$ in $L^2$ and $L^\infty$ norms.
      Bottom row: evolution of the relative errors for $v$ in $L^2$ and $L^\infty$ norms.}
      \label{f1.12}
   \end{center}
\end{figure}

\begin{figure}
   \begin{center}
      \includegraphics[scale=0.8]{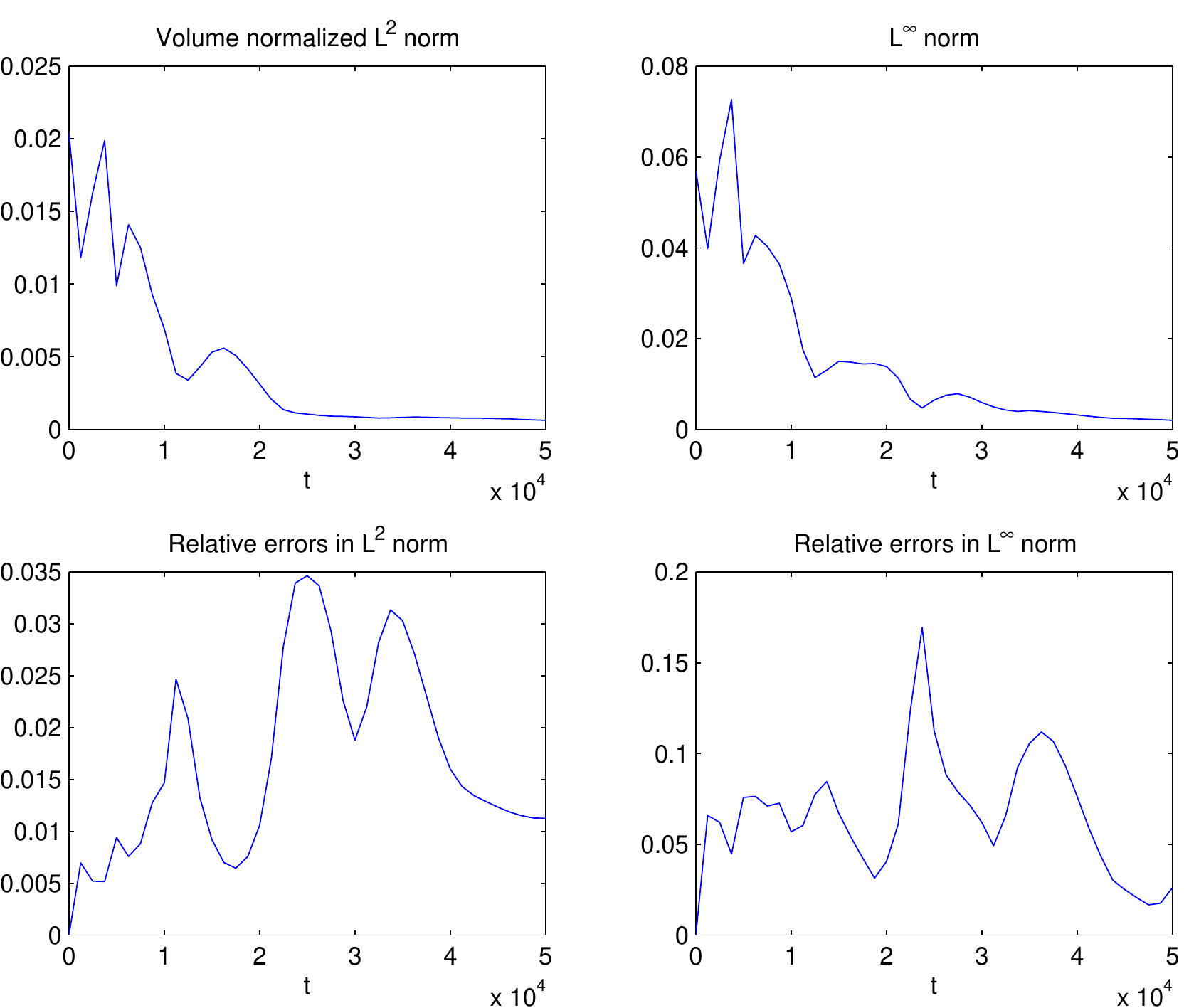}
      \caption{Top row: evolution of  the solution $w$ in $L^2$ and $L^\infty$ norms.
      Bottom row: evolution of the relative errors for $w$ in $L^2$ and $L^\infty$ norms.}
      \label{f1.13}
   \end{center}
\end{figure}

\begin{figure}
   \begin{center}
      \includegraphics[scale=0.8]{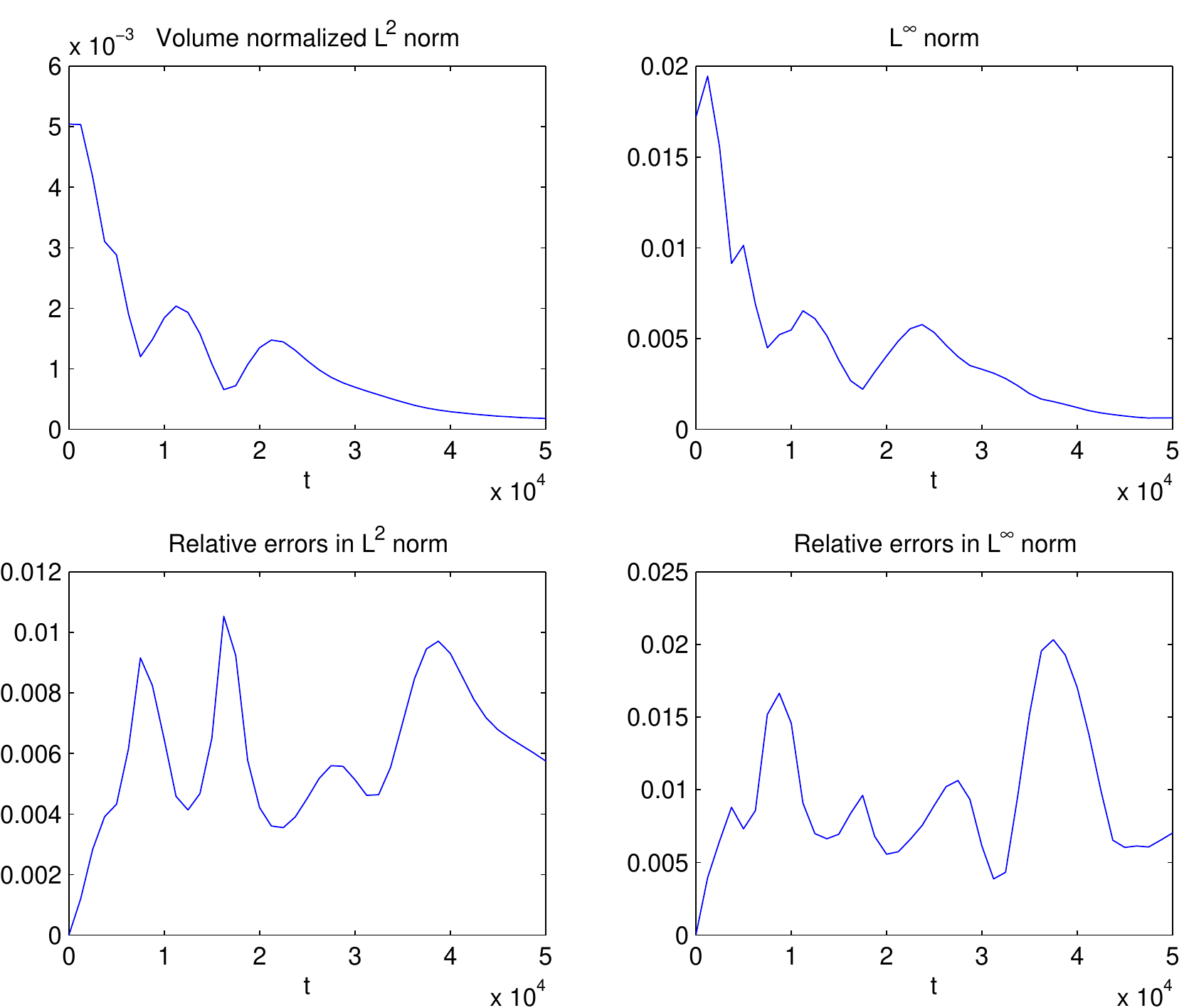}
      \caption{Top row: evolution of  the solution $\psi$ in $L^2$ and $L^\infty$ norms.
      Bottom row: evolution of the relative errors for $\psi$ in $L^2$ and $L^\infty$ norms.}
      \label{f1.14}
   \end{center}
\end{figure}

\begin{figure}
   \begin{center}
      \includegraphics[scale=0.8]{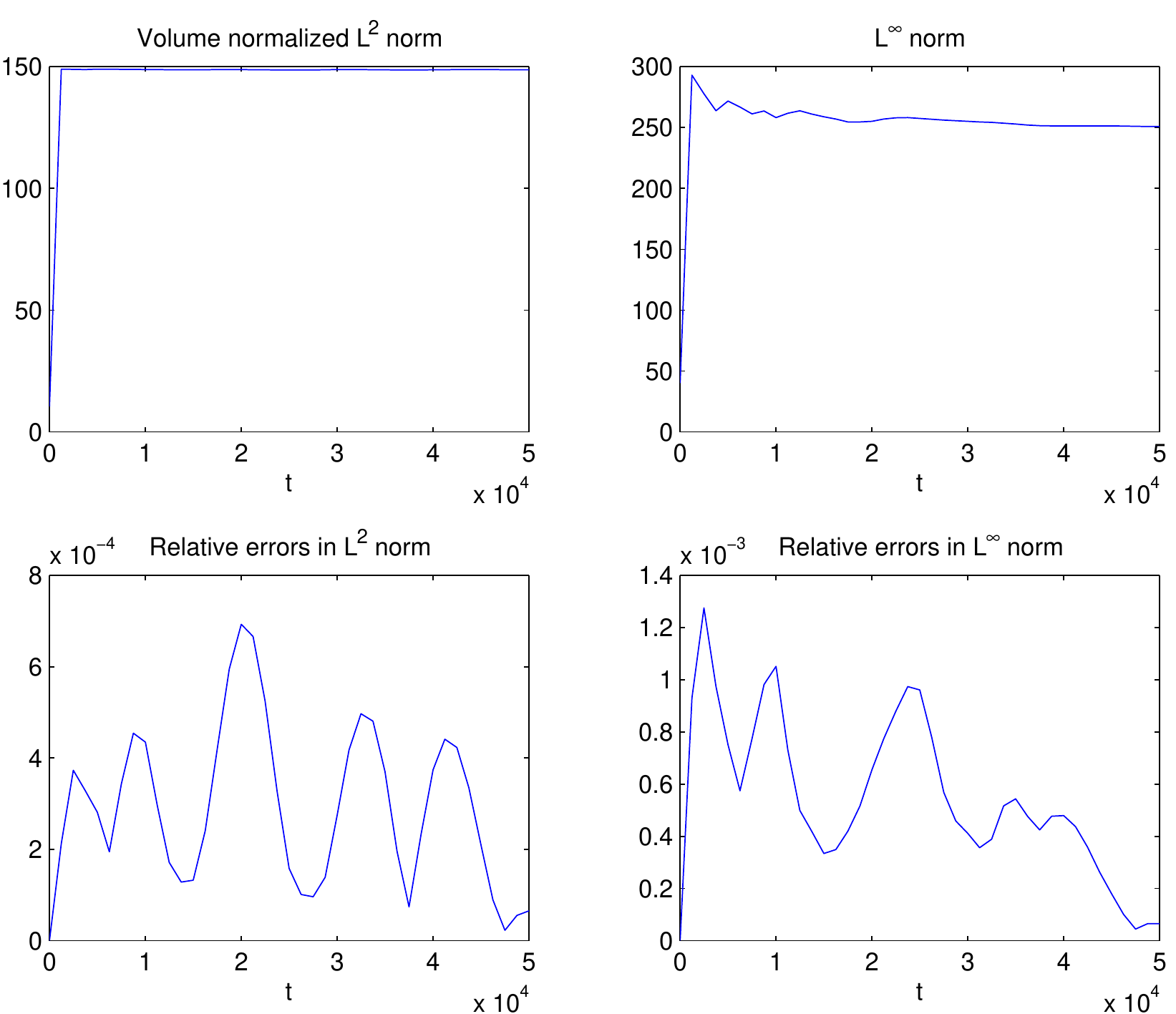}
      \caption{Top row: evolution of  the solution $\phi$ in $L^2$ and $L^\infty$ norms.
      Bottom row: evolution of the relative errors for $\phi$ in $L^2$ and $L^\infty$ norms.}
      \label{f1.15}
   \end{center}
\end{figure}

We observe that both the $L^2$ and the $L^\infty$ norms of the prognostic
variables $u$, $v$ and $\psi$ are diminishing in time. This can be explained
by the homogeneous boundary conditions imposed on the boundary of the
larger domain $\M$ and by the fact that the velocity field is flowing
out of the domain to the right, with a constant velocity $U_0$. The
relative errors of the prognostic variables, in both the $L^2$ and the
$L^\infty$ norms, are of the magnitude $O(10^{-2})$ or smaller, which means
that results on the larger domain $\M$ and the results on the middle half
domain $\M_1'$ match very well.

For the diagnostic variables, the $L^2$ and the $L^\infty$ norms of $w$
are also diminishing in time, as for the prognostic variables;
the relative errors for $w$ are large
$\left (O(10^{-1})\right )$
as compared to those for the prognostic variables, but they are still
well controlled.
The bizarre behavior of the $L^2$ and the $L^\infty$ norms of $\phi$
can be explained by the absence of an evolution equation for $\phi$ and
the lack of natural boundary conditions for
$\phi_0$ in \eqref{e2.4}. The relative errors for $\phi$ are however very
well controlled. \textbf{Note that the relative errors for $u$, $v$, and $\psi$, in both
the $L^2$ and $L^\infty$ norms, are of the order of $O(10^{-2})$, and the relative errors for $w$ are of the order of $O(10^{-1})$.}

\begin{figure}
   \begin{center}
     \epsfig{file=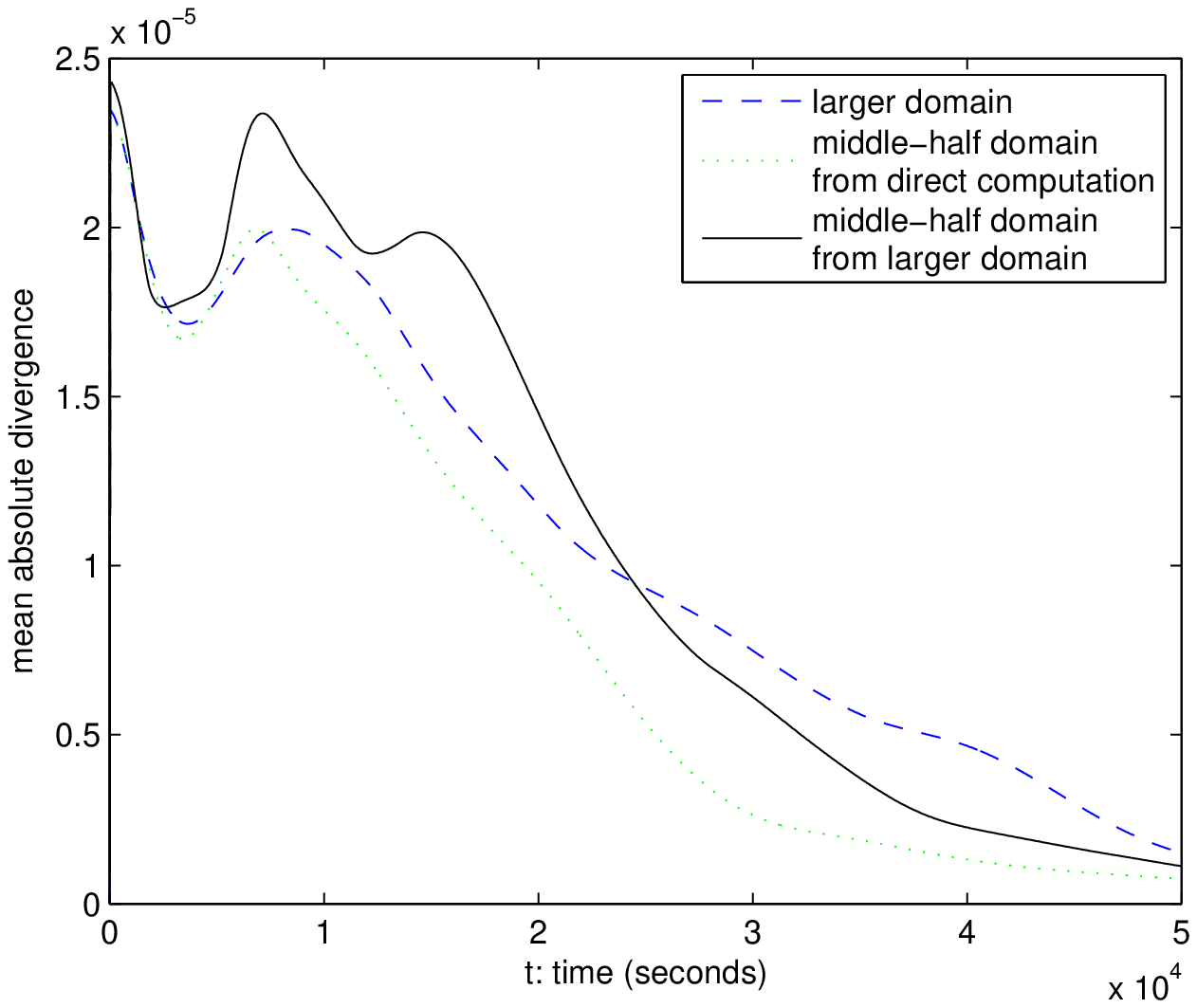,width=0.8\linewidth,clip=}
     \caption{Evolution of the mean absolute divergence for three cases.}
      \label{f1.16}
   \end{center}
\end{figure}

A graph of the absolute divergence averaged over the guest integration area is presented in Figure \ref{f1.16} for three cases: larger domain, the middle-half domain from direct computation, and the middle-half domain using the data from the larger domain. We observe that the mean absolute divergence for three cases are small and diminishing in time. This can be explained by the divergence free condition employed on the proposed numerical schemes. Furthermore, Figure \ref{f1.16} shows that the behaviors of the mean absolute divergence match well on the middle-half domain $\M_1$.

\section{Conclusions}\label{s5}
In conclusion, the absence of blowing up demonstrates that the boundary
conditions proposed in Section \ref{s2.2} are suitable for the problem,
and the numerical scheme proposed in Section \ref{s3} is stable.
The fact that the numerical results match very well on the middle half domain
$\M_1$ demonstrates the transparency property of the boundary conditions.

From the results of this idealized model it is straightforward
to outline the algorithmic path to be taken in the application of this
method to the full primitive equations. The simplest approach would be first
to re-write the  model equations so that they are formally equivalent to the
system (2.1). This involves specifying a reasonable, local mean stratification
$N^2$, and  mean zonal wind, $\bar U_0$.

Next a vertical mode decomposition is performed to identify the subcritical/supercritical mode division.
Next, the appropriate lateral boundary conditions are applied. Lastly the modal
decomposition is summed to reconstruct the boundary values of the field
variables. As necessary, the local mean stratification and zonal wind can be adjusted.

This would be superior to methods that absorb wave energy through nudging
since the artificial damping also inevitably causes the interior solution
to decay and the sponge layer, itself, induces wave reflection.

\section*{Acknowledgments}
   This work was partially supported
by the National Science Foundation under the grants NSF-DMS-0604235
and DMS-0906440, and by the Research Fund of Indiana University.

\bibliographystyle{siam}
\bibliography{biblio}

\end{document}